\documentclass{article}
\usepackage{amsmath,amssymb,url,graphicx}
\newtheorem{thm}{Theorem}[section]
\newtheorem{dfn}[thm]{Definition}
\newtheorem{prop}[thm]{Proposition}
\newtheorem{lem}[thm]{Lemma}
\newtheorem{rem}[thm]{Remark}
\numberwithin{equation}{section}
\begin{document}

\title{Cover times for sequences of reversible Markov chains on random graphs}
\author{Yoshihiro Abe \footnote{Research Institute for Mathematical Sciences, Kyoto University, Kyoto 606-8502, Japan. E-mail: yosihiro@kurims.kyoto-u.ac.jp}} 
\date{}
\maketitle

\begin{abstract}
We provide conditions that classify cover times for sequences of random walks on random graphs into two types: 
One type (Type 1) is the class of cover times that are of the order of the
maximal hitting times scaled by the logarithm of the size of vertex sets. 
The other type (Type 2) is the class of cover times that are of the order of the maximal hitting times. 
The conditions are described by some parameters determined by the underlying graphs:
the volumes, the diameters with respect to the resistance metric, the coverings or packings by balls in the resistance metric.
We apply the conditions to and classify a number of examples, such as
supercritical Galton-Watson trees,
the incipient infinite cluster of a critical Galton-Watson tree and the Sierpinski gasket graph.
\end{abstract}
Packing 
~~~~~~~~~~~~~~~~~~~~~~~~~~~~~~~~~~~~~~~~~~~~~~~~~~~~~~~~~~~~~~~Secondary 05C80
\section{Introduction and main results}
\subsection{Introduction}
\label{intro}
~~~~Let $G = (V(G), E(G))$ be a finite, connected graph and $\tau_{\text{cov}} (G)$ be the first time at which the simple random walk on $G$ visits every vertex.
The cover time for the simple random walk is defined by
$$t_{\text{cov}} (G) := \max_{x \in V(G)} E^x (\tau_{\text{cov}} (G)).$$

Cover times depend deeply on structural properties of the underlying graphs. 
Erd\H{o}s-R\'{e}nyi random graphs in several regimes are good examples. 
It is well known that as the percolation probability changes from the supercritical regime to the critical regime, 
the structure of the Erd\H{o}s-R\'{e}nyi random graph (such as the volume, the diameter) evolves.
Cooper and Frieze \cite{CF1}
and Barlow, Ding, Nachmias and Peres \cite{BDNP} estimated the cover time for the simple random walk on the 
Erd\H{o}s-R\'{e}nyi random graph in the supercritical and critical cases, respectively and
showed that the order of the cover time also evolves. 
We will investigate the relationship between cover times and structures of the underlying graphs in a more general setting.

In order to introduce our general framework, we consider the maximal hitting time defined  by
$$t_{\text{hit}} (G) := \max_{x, y \in V(G)} E^x (\tau_y (G)),$$
where $\tau_x (G)$ is the hitting time of $x$ by the simple random walk on $G.$ \\
In general, the following inequality holds for any finite, connected graphs: 
\begin{equation}
t_{\text{hit}} (G) \le t_{\text{cov}} (G) \le 2t_{\text{hit}} (G) \cdot \log |V(G)|. \label{ineq-cov-hit}
\end{equation}
The inequality on the right-hand side is often called Matthews bound (see Lemma \ref{matthews}).
In view of (\ref{ineq-cov-hit}), it is useful to classify cover times into the following two extreme types (see Definition \ref{deftype1,2} for the precise definition): \\ \\
(i) cover times that are of the order of $t_{\text{hit}} (G) \cdot \log|V(G)|$ \\
~ (we will call them Type 1), \\
(ii) cover times that are of the order of $t_{\text{hit}} (G)$ (we will call them Type 2). \\

Note that the maximal hitting time can be estimated via the volume and the diameter with respect to the resistance metric of the underlying graph
(see Lemma \ref{commute time} for the precise statement). 

In this paper, we will provide sufficient conditions that classify cover times for a sequence of random walks on random graphs into Type 1 and Type 2 in terms of
the volume, the resistance diameter and the covering or packing number of the underlying graphs 
(see section \ref{sec:1} for precise definitions of these parameters). 
We apply the conditions to many examples (see Table 1 below).
Although details of some specific cover times are already known, the novelty of this paper is that we first unify separate methods of estimating cover times into one
and add some new examples such as supercritical Galton-Watson trees and critical Galton-Watson trees conditioned to survive.

We provide intuitions for the sufficient conditions.
Roughly speaking, if one can find a packing consisting of a large number of big disjoint balls with respect to the effective resistance metric, then the cover times will be of Type 1
(Theorem \ref{type1}).
Many supercritical random graphs admit such packings. For example, we can take a family of large number of big trees as a packing
for supercritical Galton-Watson family trees and supercritical  Erd\H{o}s-R\'{e}nyi random graphs (see section \ref{sec:7}, \ref{sec:9}).

On the other hand, it can be shown that cover times will be of Type 2
if the number of balls required to cover the underlying graphs increases no more than (double) exponentially, 
as the radii of balls with respect to the resistance metric decrease exponentially
(Theorem \ref{type2}). 
A wide variety of critical random graphs and fractal graphs satisfy this property (see section \ref{sec:11}, \ref{sec:12}, \ref{sec:14}).

General bounds on cover times have been studied previously(see \cite{MCMT}, \cite{BDNP}, \cite{DLP}).
The Matthews bound (see Lemma \ref{matthews}) and the lower bound in terms of Gaussian free fields \cite{DLP} 
together with the Sudakov minoration (see Lemma \ref{sudakov}) give very useful ingredients for obtaining the condition for Type 1.
The upper bound via Gaussian free fields \cite{DLP} and the Dudley's entropy bound (see Lemma \ref{dudley}) are essential to the conditions for Type 2.

In the next subsection, we give our main results.
For a set $S$, we will write $|S|$ to denote the cardinality of $S$. 
Throughout this paper, we use $c, c^{\prime}, c_1, c_2, \dotsc $ to denote constants that does not depend on the size of $G$. 

\subsection{Main results}
\label{sec:1}
~~~~To state our main results, we first prepare some definitions.

Let $G^N = (V(G^N), E(G^N), \mu^N), N \in \mathbb{N}$  be a sequence of random weighted graphs, where $V(G^N)$ is the vertex set, $E(G^N)$ is the edge set and
$\mu^N$ is a non-negative symmetric weight function on $V(G^N) \times V(G^N)$ which satisfies $\mu^N_{xy} > 0$ if and only if $\{x,y\} $ $\in E(G^N).$
We assume that these weighted graphs are defined on a common probability space with
a probability measure $\textbf{P}$ and that $G^N$ is a finite, connected graph, $\textbf{P}$-a.s. 
In this paper, the following four parameters (volume, resistance diameter, packing number, covering number) play important roles in estimating cover times. \\
The volume of $G^N$ is defined by 
$$\mu^N (G^N) := \sum_{x, y \in V(G^N)} \mu_{xy}^N.$$ 
The effective resistance is a powerful tool for studying random walks on weighted graphs (see Lemma \ref{commute time}).
For $x,y \in V(G^N), x \neq y,$ we define the effective resistance between $x$ and $y$ by
$$R_{\text{eff}}^N (x,y)^{-1} := \inf \{\mathcal{E}^N(f,f) : f \in \mathbb{R}^{V(G^N)}, f(x) = 1, f(y) = 0 \},$$
where 
$\mathcal{E}^N (f,g) := \frac{1}{2} \displaystyle \sum_{\begin{subarray}{c} u,v \in V(G^N) \\ \{u,v\} \in E(G^N) \end{subarray} } \mu^N_{uv} (f(u) - f(v))(g(u) - g(v)), f,g \in 
\mathbb{R}^{V(G^N)}.$ \\
If we define $R_{\text{eff}}^N (x,x) = 0$ for all $x \in V(G^N),$ it is known that $R_{\text{eff}}^N (\cdot, \cdot)$ is a metric on $V(G^N).$
The resistance diameter is defined by 
$$\text{diam}_{R}(G^N) := \displaystyle \max_{x,y \in V(G^N)} R^N_{\text{eff}}(x,y).$$
We define the resistance ball with radius $r$ centered at $x \in V(G^N)$ by 
$$B_{\text{eff}}^N (x, r) := \{y \in V(G^N) : R_{\text{eff}}^N (x, y) \le r\}.$$
We call a family of resistance balls $\{B_{\text{eff}}^N (x_1, r_1), \cdots, B_{\text{eff}}^N (x_{m}, r_m) \}$ a packing for $G^N$ if these resistance balls are disjoint with each other. \\
The packing number for $(G^N, r)$ is defined by
\begin{align*}
n_{\text{pac}} (G^N, r) := \max \Big \{m &\ge 1 : ~\text{there exist} ~x_1, \cdots, x_m \in V(G^N) ~\text{such that} \\
                                          &\{B_{\text{eff}}^N (x_1, r), \cdots, B_{\text{eff}}^N (x_m, r) \}~ \text{is a packing for}~ G^N \Big \}.
\end{align*}
We call a family of resistance balls $\{B_{\text{eff}}^N (x_1, r_1), \cdots, B_{\text{eff}}^N (x_{m}, r_m) \}$ a covering for $G^N$ if 
$$V(G^N) \subset \displaystyle \bigcup_{k = 1}^m B_{\text{eff}}^N (x_k, r_k).$$
The covering number for $(G^N, r)$ is defined by
\begin{align*}
n_{\text{cov}} (G^N, r) := \min \Big \{m &\ge 1 : ~\text{there exist} ~x_1, \cdots, x_m \in V(G^N) ~\text{such that} \\
                                                       &\{B_{\text{eff}}^N (x_1, r), \cdots, B_{\text{eff}}^N (x_m, r) \}~ \text{is a covering for}~ G^N \Big \}.
\end{align*}

The discrete time random walk on $G^N$ is the Markov chain $((X_n)_{n \geqslant 0}, P^x, x \in V(G^N))$ with transition probabilities $(p(x,y))_{x,y \in 
V(G^N)}$
defined by $p(x,y)$ $:= \mu_{xy}^N/\mu_x^N$, where $\mu_x^N := \sum_{y \in V(G^N)} \mu_{xy}^N$.
Let $\tau_{\text{cov}} (G^N)$ be the first time at which the random walk visits every vertex of $V(G^N)$.
We define the cover time for the random walk on $G^N$ as follows:
$$t_{\text{cov}}(G^N) := \max_{x \in V(G^N)} E^x(\tau_{\text{cov}}(G^N)).$$
We also define the maximal hitting time for the random walk on $G^N$ by
$$t_{\text{hit}} (G^N) := \max_{x, y \in V(G^N)} E^x(\tau_y (G^N)),$$
where $\tau_x (G^N)$ is the hitting time of $x \in V(G^N)$ by the random walk on $G^N$.
We give the precise definitions of types for a sequence of cover times.
\begin{dfn} \label{deftype1,2} 
(1) A sequence of cover times $(t_{\text{cov}} (G^N) )_{N \in \mathbb{N} }$ is Type 1 if
\begin{equation}
\lim_{\lambda \to \infty} \liminf_{N \to \infty} \textbf{P} \left( \lambda^{-1} \le \frac{t_{\text{cov}} (G^N)}{t_{\text{hit}} (G^N) \cdot \log |V(G^N)|} \le 2 \right) = 1. \label{def1}
\end{equation}
(2) A sequence of cover times $(t_{\text{cov}} (G^N) )_{N \in \mathbb{N} }$ is Type 2 if 
\begin{equation}
\lim_{\lambda \to \infty} \liminf_{N \to \infty} \textbf{P} \left( 1 \le \frac{t_{\text{cov}} (G^N)}{t_{\text{hit}} (G^N)} \le \lambda \right) = 1. \label{def2}
\end{equation}
\end{dfn}
\begin{rem} \label{remgeneralbound}
By (\ref{ineq-cov-hit}),
the upper bound of the event in (\ref{def1}) and the lower bound of the event in (\ref{def2}) always hold. 
\end{rem}

We are now ready to state our main theorems.
We first state the sufficient condition for cover times to be Type 1.
We will say that a sequence of events $(B_N)_{N \ge 0} $ holds with high probability (abbreviated to w.h.p.) if $\lim_{N \to \infty} \textbf{P}(B_N) = 1$. 
\begin{thm}  \label{type1}
(1) Suppose there exist $c_1, c_2> 0$ and functions $v,r : \mathbb{N} \to [0,\infty) $ with $ \lim_{N \to \infty} v(N) = \infty$ such that w.h.p., the following holds: 
\begin{equation}
\log |V(G^N)| \le c_1 \log v(N),~\text{diam}_R (G^N) \le c_2 r(N). \label{condition-of-type1-0}
\end{equation}
Then there exists $c_3 >0$ such that w.h.p.,
$$t_{\text{cov}}(G^N)/\mu^N(G^N) \le c_3 r(N)\log v(N).$$

(2) Suppose that there exist $c_4, c_5 > 0$ and functions $v,r : \mathbb{N} \to [0,\infty) $ with $ \lim_{N \to \infty} v(N) = \infty$ such that w.h.p.,
\begin{equation}
\log \{n_{\text{pac}} (G^N, c_4 r(N)) \} \ge c_5 \log v(N). \label{condition-of-type1}
\end{equation}
Then there exists $c_6 > 0$ such that w.h.p.,
$$t_{\text{cov}}(G^N)/\mu^N(G^N) \ge c_6 r(N)\log v(N).$$

(3) Under conditions (\ref{condition-of-type1-0}) and (\ref{condition-of-type1}), $(t_{\text{cov}} (G^N))_{N \in \mathbb{N}}$ is Type 1. \\
\end{thm}

We next state sufficient conditions for cover times to be Type 2. 
\begin{thm} \label{type2}
(1) Suppose that there exist functions $v,r : \mathbb{N} \to [0,\infty) $ with $ \lim_{N \to \infty} v(N) = \infty$ and a function $p : [1, \infty) \to [0, 1]$ with
$\lim_{\lambda \to \infty} p(\lambda) = 0$ satisfying the following for all $\lambda \ge 1$ and sufficiently large $N \in \mathbb{N}:$
\begin{equation}
\textbf{P}(\mu^N(G^N) \le \lambda v(N)) \ge 1 - p(\lambda), \label{condition-of-type2-(1)-1} 
\end{equation}
and there exists a random non-increasing sequence $(\ell_k^N)_{k \ge 0}$ satisfying
$\ell_0^N = \text{diam}_R(G^N), \ell_{k_0^N - 1}^N > 0$ and $\ell_{k_0^N}^N = 0$ for some $k_0^N \in \mathbb{N}$ such that
\begin{equation} 
\textbf{P}\Big (\sum_{k = 1}^{k_0^N} \sqrt{\ell_{k-1}^N \log \{ n_{\text{cov}} (G^N, \ell_k^N) \}}  \le \lambda \sqrt{r(N)} \Big)  
\ge 1 - p(\lambda). \label{condition-of-type2-(1)-2}
\end{equation}
Then there exists $c>0$ such that for all $\lambda \ge c$ and sufficiently large $N \in \mathbb{N}$,
\begin{equation} \label{cov-estimate-type2}
\textbf{P}(t_{\text{cov} }(G^N) > \lambda v(N) r(N)) \le \inf_{0 < \theta < 1} \left \{ p((\lambda/c)^{\theta}) + p \left((\lambda/c)^{\frac{1 - \theta}{2} } \right) \right \}.
\end{equation}
(2) Suppose that there exist functions $v,r : \mathbb{N} \to [0,\infty) $ with $ \lim_{N \to \infty} v(N) = \infty$ and a function $p : [1, \infty) \to [0, 1]$ with
$\lim_{\lambda \to \infty} p(\lambda) = 0$ satisfying the following for all $\lambda \ge 1$ and sufficiently large $N \in \mathbb{N}:$
\begin{equation}
\textbf{P}(\mu^N(G^N) < {\lambda}^{-1} v(N)) \le p(\lambda), ~\textbf{P}(\text{diam}_R(G^N) < {\lambda}^{-1} r(N) ) \le p(\lambda). \label{condition-of-type2-(3)}
\end{equation}
Then there exists $c > 0$ such that for all $\lambda \ge c$ and sufficiently large $N \in \mathbb{N}$,
$$\textbf{P}(t_{\text{cov} }(G^N) < {\lambda}^{-1} v(N) r(N)) \le \inf_{0 < \theta < 1} \Big \{ p \Big (\Big (\frac{\lambda}{c} \Big)^{\theta} \Big) +
p \Big ( \Big (\frac{\lambda}{c} \Big)^{1 - \theta} \Big) \Big \}.$$
(3) Under the conditions (\ref{condition-of-type2-(1)-1}), (\ref{condition-of-type2-(1)-2}) and (\ref{condition-of-type2-(3)}), 
$(t_{\text{cov}} (G^N))_{N \in \mathbb{N}}$ is Type 2. \\
\end{thm}
\begin{rem} \label{type2rems}
(1) In general, we can not replace (\ref{cov-estimate-type2}) by the statement that $t_{\text{cov} }(G^N) \le cv(N)r(N)$ w.h.p., for some $c > 0$ (see Proposition \ref{cgwprop2}).
We thus state Theorem \ref{type1} and Theorem \ref{type2} in a slightly different way. \\
(2) If the conditions (\ref{condition-of-type1-0}) and (\ref{condition-of-type1}) in Theorem \ref{type1} hold $\textbf{P}$-almost surely for sufficiently large $N \in \mathbb{N},$ 
the results of Theorem \ref{type1} also hold $\textbf{P}$-almost surely for sufficiently large $N \in \mathbb{N}.$ \\
(3) If the events of (\ref{condition-of-type2-(1)-1}), (\ref{condition-of-type2-(1)-2}) and (\ref{condition-of-type2-(3)}) 
in Theorem \ref{type2} hold $\textbf{P}$-almost surely for sufficiently large $N \in \mathbb{N},$
the results of Theorem \ref{type2} also hold $\textbf{P}$-almost surely for sufficiently large $N \in \mathbb{N}$ ($\lambda$ will be replaced by some constants). \\
(4) On some class of planar graphs, the condition (\ref{condition-of-type1}) always holds;
Let $(G^N)_{N \ge 0}$ be a sequence of $\textbf{P}$-a.s. finite, planar connected random graphs with maximum degree $c > 0$ and $\mu_{xy}^N = 1$ for all $\{x, y \} \in E(G^N).$
Suppose that there exists $c_7 > 0$ and a function $v : \mathbb{N} \to [0,\infty) $ with $ \lim_{N \to \infty} v(N) = \infty$ such that w.h.p., 
$\log |V(G^N)| \ge c_7 \log v(N).$
Then by Lemma 3.1 of \cite{JS}, (\ref{condition-of-type1}) holds with the function $v$ and $r(N) = \log v(N).$ \\
(5) Typically, we take an exponentially decreasing sequence as $(\ell_k^N)_{k \ge 0}$ in (\ref{condition-of-type2-(1)-2}) (for example, $\ell_k^N = \frac{\text{diam}_R (G^N)}{2^k}$). 
\end{rem} 
Applying these theorems, we will estimate and classify cover times for several specific random graphs. We give a list of the results in Table 1. 
We explain the notation in Table 1.
The notation $m$ is the mean of the offspring distribution of the corresponding branching process. 
Supercritical Erd\H{o}s-R\'{e}nyi random graphs I, II have the percolation probability $c/N$, $f(N)/N$
respectively, where $c > 1$ is a constant and $\lim_{N \to \infty} f(N)/ \sqrt{N} = \lim_{N \to \infty} \log N/ f(N) = 0.$  
`IIC' is the abbreviation of `incipient infinite cluster' and $p_N$ is the survival probability up to $N$ level (see subsection \ref{sec:11}).

\begin{center}
Table 1: Orders of volumes and cover times for random graphs and types of cover times
\end{center}
\begin{tabular}{|l|c|c|c|} \hline
Random graph & Volume & Cover time & Type \\ \hline \hline
Supercritical Galton-Watson family trees & $m^N$ & $N^2 m^N$ & 1 \\ \hline
Supercritical Erd\H{o}s-R\'{e}nyi random graphs I & $N$ & $N (\log N)^2$ & 1 \\ \hline 
Supercritical Erd\H{o}s-R\'{e}nyi random graphs II & $N f(N)$ & $N \log N$ & 1 \\ \hline 
The IIC for critical Galton-Watson family tree & $N p_N^{- 1}$ & $N^2 p_N^{- 1}$ & 2 \\ \hline 
Critical Erd\H{o}s-R\'{e}nyi random graphs & $N^{2/3}$ & $N$ & 2  \\ \hline
The range of random walk in $\mathbb{Z}^d, d \ge 5$ & $N$ & $N^2$ & 2  \\ \hline
Sierpinski gasket graphs & $3^N$ & $5^N$ & 2  \\ \hline
\end{tabular} \\ \\

Concerning the IIC for critical Galton-Watson family trees,
Aldous \cite{Al} and Barlow, Ding, Nachmias and Peres \cite{BDNP} have estimated
the cover times for critical Galton-Watson family trees for finite variance offspring distributions. 
Our result extends these results to the case where the offspring distribution is in the domain of attraction of a stable law with index $\alpha \in (1, 2].$ 
Our result clarifies that the cover time for the IIC depends on the survival probability of the branching process up to some level.

In addition to this example, we give new estimates on cover times for supercritical Galton-Watson family trees, the range of random walk in $\mathbb{Z}^d, d \ge 5$
and Sierpinski gasket graphs.

Note that for supercritical Erd\H{o}s-R\'{e}nyi random graphs, better estimates are already known \cite{CF1, J}
and that for critical Erd\H{o}s-R\'{e}nyi random graphs, the correct order is already known \cite{BDNP}.
We cite these examples to compare Type 1 and Type 2.

In Section \ref{sec:8}, we will estimate the cover time for the largest supercritical percolation cluster inside a box in $\mathbb{Z}^d, d \ge 2.$ However, we are not able to obtain 
the correct order (see Remark \ref{remscpc}).

Note that there are graphs where the cover times can not be classified as either Type 1 or Type 2.
For example, let $G^N$ be a deterministic graph with unit weights
consisting of a complete graph with $N$ vertices and $a_N$ other vertices, each attached by a single edge to a distinct vertex of the complete graph, 
where $a_N$ is a positive number satisfying $2 \le a_N \le N$.
One can show that $\text{diam}_R (G^N) = 2 + 2/N,$ $n_{\text{pac}} (G^N, \ell) \ge a_N$ for all $0 \le \ell \le 1$ 
and $n_{\text{cov}} (G^N, \text{diam}_R (G^N)/2^k) \le a_N + 1$ for all $1 \le k \le \lfloor \log_2 N \rfloor.$ 
By Theorem \ref{type1} (2), Lemma \ref{commute time} and Lemma \ref{barlow-ding-nachmias-peres} below, we have for some $c, c^{\prime} > 0,$
$$c \cdot t_{\text{hit}} (G^N) \cdot \log a_N \le t_{\text{cov}} (G^N) \le c^{\prime} \cdot t_{\text{hit}} (G^N) \cdot \log a_N.$$
This implies that if $\displaystyle \lim_{N \to \infty} a_N = \infty$ and $\displaystyle \lim_{N \to \infty} \frac{\log a_N}{\log N} = 0$, then the sequence of cover times
$(t_{\text{cov}} (G^N))_{N \in \mathbb{N}}$ is neither Type 1 nor Type 2.

We give the outline of this paper. 
In Section \ref{sec:2}, we prove Theorem \ref{type1} and Theorem \ref{type2}. 
In Section \ref{sec:6}, using Theorem \ref{type1} and Theorem \ref{type2}, we estimate and classify cover times for the examples in Table 1. 

\section{Proof of Theorem \ref{type1} and Theorem \ref{type2}} 
\label{sec:2}
~~~~In this section, we prove Theorem \ref{type1} and Theorem \ref{type2}.
\subsection{Known results}
\label{sec:3}
~~~~We state some known results on cover times and Gaussian free fields that we will use in this paper. 

Throughout the following lemmas, $G = (V(G), E(G))$ will be a finite, connected graph and $\mu$ will be the weight function with $\mu (G) := \sum_{x, y \in V(G)} \mu_{xy}.$
Let $\{\eta_x \}_{x \in V(G)}$ be the Gaussian free field on $G$ defined on a probability space with a probability measure $\mathbb{P}.$ \\
Recently, Ding, Lee and Peres \cite{DLP} proved the following surprising result, which says that cover times have a close relationship with Gaussian free fields.
\begin{lem}(\cite{DLP}, Theorem 1.9 and Theorem (MM))
\label{ding-lee-peres}
There exist $c_1, c_2 > 0$ such that
$$c_1 \cdot \mu(G) \cdot \Bigg(\mathbb{E} \max_{x \in V(G) } \eta_x \Bigg)^2 \le t_{\text{cov} }(G) \le c_2 \cdot \mu(G) \cdot \Bigg(\mathbb{E} \max_{x \in V(G) } \eta_x \Bigg)^2.$$
\end{lem}

The following commute time identity is well-known and useful for estimating the maximal hitting time. See, for instance, Theorem 2.1 of \cite{CRRST}
or Proposition 10.6 of \cite{MCMT}.

\begin{lem} \label{commute time}
Let $\tau_x$ be the hitting time of $x \in V(G)$ by the random walk on $G$. For all $x, y \in V(G)$,
$$E^x (\tau_y) + E^y (\tau_x) = \mu(G) R_{\text{eff}} (x, y).$$
In particular,
$$\frac{1}{2} \mu (G) \text{diam}_R(G) \le t_{\text{hit}} (G) \le \mu (G) \text{diam}_R (G).$$
\end{lem}

Fix $x, y \in V(G)$. $\Pi$ is an edge-cutset between $x$ and $y$ if $\Pi$ is a subset of $E(G)$ such that every path from $x$ to $y$
has an edge belonging to $\Pi$. 
The following Nash-Williams inequality is useful for obtaining lower bounds on effective resistances. See, for example, Proposition 9.15 of \cite{MCMT}.
\begin{lem} \label{Nash-Williams} Fix $x, y \in V(G)$.
Let $(\Pi_k)_{k \ge 1}$ be a sequence of edge-cutsets between $x$ and $y$ with $\Pi_k \cap \Pi_{\ell} = \emptyset$ for all $k \neq \ell$. Then,
$$R_{\text{eff}} (x, y) \ge \sum_{k \ge 1} \big(\sum_{\{u, v \} \in \Pi_k} \mu_{uv} \big)^{- 1}.$$
\end{lem}

\subsection{Proof of Theorem \ref{type1}}
\label{sec:4}
~~~We provide the proof of Theorem \ref{type1}. The following lemma is known as the Matthews bound. See, for example, Theorem 11.2 of \cite{MCMT}
(see also the original work of Matthews \cite{Ma}).
\begin{lem}  \label{matthews}
Let $(X_n)_{n \ge 0}$ be an irreducible Markov chain on a finite state space $V$ and $t_{\text{cov}}, t_{\text{hit}}$ be its cover time and maximal hitting time, respectively. 
Then, 
$$t_{\text{cov}} \le t_{hit} \cdot (\log |V| + 1).$$
\end{lem}

We also use the next fact, called Sudakov minoration. See, for instance, Lemma 2.1.2 of \cite{TGC}.
\begin{lem}
\label{sudakov}
Let $\{\eta_x \}_{x \in V(G)}$ be a Gaussian free field on a weighted graph $G$.
There exists $c > 0$ such that for all $V^{'} \subset V(G)$, 
$$\mathbb{E} \max_{x \in V^{'}} \eta_x \ge c \Big(\min_{\begin{subarray}{c} y, z \in V^{'} \\ y \neq z \end{subarray} } \sqrt{R_{\text{eff}} (y, z)}  ~\Big) \sqrt{\log |V^{'}|}.$$
\end{lem}
\textit{Proof of Theorem \ref{type1}.}
We first prove (1).
By Lemma \ref{commute time} and (\ref{condition-of-type1-0}), we get w.h.p.,
\begin{equation}
t_{\text{hit}} (G^N) \le \mu^N (G^N) \cdot \text{diam}_R (G^N) \le c_2 \mu^N(G^N) r(N). \label{ineq-of-pf-of-thm1.1-1} 
\end{equation}
So, using Lemma \ref{matthews}, (\ref{condition-of-type1-0}) and (\ref{ineq-of-pf-of-thm1.1-1}) , we have that w.h.p.,
\begin{align*}
t_{\text{cov}} (G^N) &\le t_{\text{hit}} (G^N) \cdot (\log |V(G^N)| + 1)  \\
                          & \le 2 c_1 c_2 \mu^N (G^N) r(N) \log v(N).   
\end{align*}
Next, we prove (2).
Let $x_1, \cdots, x_{n_{\text{pac}} (G^N, c_4 r(N))}$ be vertices satisfying that
the set of resistance balls $\{B_{\text{eff}}^N (x_k, c_4 r(N)) : 1 \le k \le n_{\text{pac}} (G^N, c_4 r(N)) \}$ is a packing for $G^N.$
Set $V^{'} := \{ x_1, \cdots, x_{n_{\text{pac}} (G^N, c_4 r(N))} \}.$
Using (\ref{condition-of-type1}), Lemma \ref{ding-lee-peres} and Lemma \ref{sudakov}, we have that there exist $c_7, c_8 > 0$ such that w.h.p., 
\begin{align}
t_{\text{cov}} (G^N) &\ge c_7 \mu^N (G^N) \Big (c_8  \sqrt{c_4 r(N)} \sqrt{\log \{n_{\text{pac}} (G^N, c_4 r(N)) \}} \Big )^2  \notag \\
                           &\ge c_4 c_5 c_7 c_8^2 ~\mu^N (G^N) r(N) \log v(N).  \label{ineq-of-pf-of-thm1.1-2}                                                
\end{align}
The inequalities (\ref{condition-of-type1-0}), (\ref{ineq-of-pf-of-thm1.1-1}) and (\ref{ineq-of-pf-of-thm1.1-2}) imply the conclusion of (3). 
$\Box$ 

\subsection{Proof of Theorem \ref{type2}}
\label{sec:5}
~~~~We prove Theorem \ref{type2}. 
The following fact is a minor extension of Theorem 1.1 of \cite{BDNP} and provides useful general upper bounds on cover times.
\begin{lem}
\label{barlow-ding-nachmias-peres}
Let $G = (V(G), E(G))$ be a graph and $\mu$ be the weight function with $\mu (G) := \sum_{x, y \in V(G)} \mu_{xy}.$ \\
Let $(\ell_k)_{k \ge 0}$ be a non-increasing sequence with $\ell_0 = \text{diam}_R (G), \ell_{k_0 - 1} > 0$ and $\ell_{k_0} = 0$ for some $k_0 \in \mathbb{N}.$ \\
There exists $c > 0$ such that
$$t_{\text{cov} }(G) \le c \bigg(\sum_{k = 1}^{k_0} \sqrt{\ell_{k-1} \log \{n_{\text{cov}} (G, \ell_k) \} } \bigg)^2 \cdot \mu(G).$$
\end{lem}
Lemma \ref{barlow-ding-nachmias-peres} follows from the following result. See, for example, Theorem 11.17 of \cite{LT}.
\begin{lem} \label{dudley}
Let $I$ be a finite set and $\{\eta_x \}_{x \in I}$ be a Gaussian process. 
Set $d (x,y) := \sqrt{\mathbb{E} (\eta_x - \eta_y)^2}$ and
\begin{align*}
n(I, d, \ell) := \min \{m \ge 1 :  ~&\text{there exist}~x_1, \cdots, x_m \in I \\
&\text{such that}~I \subset \bigcup_{k = 1}^m \{y \in I : d(x_k, y) \le \ell \} \}.
\end{align*}
Then there exists $c > 0$ such that 
$$\mathbb{E} \max_{x \in I} \eta_x \le c \int_0^{\infty} \sqrt{\log \{n (I, d, \ell) \} } d\ell.$$
\end{lem}
\textit{Proof of Lemma \ref{barlow-ding-nachmias-peres}.}
Let $\{\eta_x \}_{x \in V(G)}$ be a Gaussian free field on $G.$ 
Note that $d(x,y) = \sqrt{\mathbb{E} (\eta_x - \eta_y)^2} = \sqrt{R_{\text{eff}} (x, y)}.$ \\
In particular, $n(V(G), d, \ell) = n_{\text{cov}} (G, \ell^2).$  \\
Since $n_{\text{cov}} (G, \ell)$ is non-increasing with respect to $\ell,$ we have
\begin{align}
&~~~\int_0^{\infty} \sqrt{\log \{n (V(G), d, \ell) \} } d\ell \notag \\
&\le \int_0^{\infty} \sqrt{\log \{n_{\text{cov}} (G, \ell^2) \} } d\ell \notag \\
&\le \sum_{k = 1}^{k_0} \int_{\sqrt{\ell_k}}^{\sqrt{\ell_{k-1}}} \sqrt{\log \{n_{\text{cov}} (G, \ell^2) \} } d\ell \notag \\
&\le \sum_{k = 1}^{k_0} \sqrt{\ell_{k-1} \log \{n_{\text{cov}} (G, \ell_k) \} }. \label{dudley-bound}
\end{align}
Lemma \ref{ding-lee-peres}, Lemma \ref{dudley} and  (\ref{dudley-bound}) imply the conclusion. $\Box$

\textit{Proof of Theorem \ref{type2}.}
First, we prove (1). Fix $\lambda \ge 1,$ sufficiently large $N \in \mathbb{N}$ and $\theta \in (0, 1).$
Set 
\begin{equation*}
B := \Big \{ \sum_{k = 1}^{k_0^N} \sqrt{\ell_{k-1}^N \log \{ n_{\text{cov}} (G^N, \ell_k^N) \}}  \le \lambda^{\frac{1-\theta}{2}} \sqrt{r(N)} \Big \}.
\end{equation*}
By (\ref{condition-of-type2-(1)-1}), (\ref{condition-of-type2-(1)-2}) and Lemma \ref{barlow-ding-nachmias-peres}, we have for some $c_1 > 0$ that
\begin{align*}
&~~~\textbf{P} (t_{\text{cov}} (G^N) > c_1 \lambda v(N) r(N) ) \\
&\le \textbf{P} (\mu^N(G^N) > \lambda^{\theta} v(N) ) + \textbf{P} (B^c)\\
&\le p(\lambda^{\theta}) + p(\lambda^{\frac{1 - \theta}{2}}),
\end{align*}
which implies the conclusion of (1).

Next, we prove (2). Fix $\lambda \ge 1,$ sufficiently large $N \in \mathbb{N},$ and $\theta \in (0, 1).$
By (\ref{condition-of-type2-(3)}), Lemma \ref{commute time} and the fact that $t_{\text{cov}} (G^N) \ge t_{\text{hit}} (G^N)$ $\textbf{P}$-a.s., we have that
\begin{align*}
&~~~\textbf{P} \left(t_{\text{cov}} (G^N) < \frac{\lambda^{- 1}}{2} v(N) r(N) \right ) \\
&\le \textbf{P} (\mu^N (G^N) < \lambda^{- \theta} v(N) ) + \textbf{P} (\text{diam}_R(G^N) < \lambda^{- (1 - \theta )} r(N) ) \\
&\le p(\lambda^{\theta}) + p(\lambda^{1 - \theta}),
\end{align*}
which implies the conclusion of (2). 

Using Lemma \ref{commute time} and the results of (1) and (2), we can easily obtain the conclusion (3). We omit the detail. $\Box$ 
\section{Examples}
\label{sec:6}
~~~In this section, we estimate and classify cover times for a number of specific random graphs by using Theorem \ref{type1} and Theorem \ref{type2}. 
Given a graph $G$, we will write $d_G (x, y)$ to denote the graph distance between $x$ and $y$ in the graph $G$.  
From Subsection 3.1 to 3.7, we assume that $\mu_{xy}^N = 1$ for all $\{x, y \} \in E(G^N)$ and $N \in \mathbb{N}$ $\textbf{P}$-a.s.
\subsection{Supercritical Galton-Watson family trees}
\label{sec:7}
~~~Let $(Z_N)_{N \ge 0}$ be a Galton-Watson process defined on a probability space with probability measure $\mathbb{P}$ and $\mathcal{T}$ be its family tree. 
We assume that $m := \mathbb{E} (Z_1) \in (1, \infty).$
$\mathcal{T}_{\le N}$ and $\mathcal{T}_N$ are the first $N$ generations and the set of $N$-th generation of $\mathcal{T}$ respectively. 
In particular, $Z_N = |\mathcal{T}_N|.$
$\Tilde{\mathcal{T}}_N$ is a set of vertices among $N$-th generation that have infinite line of descent. 
We consider the conditional measure
$\textbf{P} := \mathbb{P}(~ \cdot ~ |~Z_n \neq 0 ~\text{for all}~ n \in \mathbb{N})$. 
We prove the following proposition.
\begin{prop} \label{scgwprop} There exist $c_1, c_2 > 0$ such that $\textbf{P}$-a.s., for sufficiently large $N \in \mathbb{N},$ 
$$c_1 N^2 \le t_{cov}(\mathcal{T}_{\le N}) / |E(\mathcal{T}_{\le N})| \le c_2 N^2,$$
and $(t_{cov}(\mathcal{T}_{\le N}))_{N \in \mathbb{N}}$ is Type 1.
\end{prop} 
In the proof, we use the following well-known fact. See, for example, Theorem 1 (page 49), Theorem 3 (page 30) and Lemma 4 (page 31) of \cite{AN}.
\begin{lem} \label{scgwlem}
Let $(Z_N)_{N \ge 0}$ be a Galton-Watson process with mean $m \in (1, \infty).$ \\
(1) Set $\tilde{Z}_N := |\tilde{\mathcal{T}}_N|.$ Under the probability measure $\mathbb{P} (\cdot | Z_n \neq 0~ \text{for all}~ n \in \mathbb{N}),$
$(\tilde{Z}_N)_{N \ge 0}$ is a Galton-Watson process whose offspring distribution has generating function 
$$\tilde{f} (s) = \frac{f((1-q)s + q) - q}{1 - q},$$
where $f$ is the generating function of $Z_1$ and $q := \mathbb{P} (Z_n = 0~ \text{for some}~ n \in \mathbb{N}).$ \\
(2) There exist a sequence of constants $(C_N)_{N \in \mathbb{N}}$ with $\displaystyle \lim_{N \to \infty} C_N = \infty$ 
and $\displaystyle \lim_{N \to \infty} \frac{C_{N + 1}}{C_N} = m$ and a random variable $W$ such that
$$\lim_{N \to \infty} \frac{Z_N}{C_N} = W ~\mathbb{P} \text{-a.s.}, ~\mathbb{P} (W < \infty) = 1 ~\text{and}~\mathbb{P} (W = 0) = q.$$
\end{lem}

\textit{Proof of Proposition \ref{scgwprop}.} 
We check almost-sure versions of (\ref{condition-of-type1-0}) and (\ref{condition-of-type1}) in Theorem \ref{type1} with $\log v(N) = r(N) = N.$ \\
By the Chebyshev inequality, we have for all $\alpha > m,$ 
$$\textbf{P} (|\mathcal{T}_{\le N}| > \alpha^N) \le \frac{\textbf{E} (|\mathcal{T}_{\le N}|)}{\alpha^N} \le \frac{1}{1-q} \cdot \frac{m}{m - 1} \cdot \Big (\frac{m}{\alpha} \Big )^N.$$
So, by the Borel-Cantelli lemma, $|\mathcal{T}_{\le N}| \le \alpha^N$ for sufficiently large $N \in \mathbb{N},$ $\textbf{P}$-a.s.
Since $R_{\text{eff}}^N (x,y) = d_{\mathcal{T}_{\le N}} (x,y)$ for all $x,y \in \mathcal{T}_{\le N},$ we get $\text{diam}_R (\mathcal{T}_{\le N} ) \le 2N,$
$\textbf{P}$-a.s.
We set $V^{'} := \{g_N(v) : v \in \Tilde{\mathcal{T}}_{\lfloor \frac{N}{2} \rfloor} \}$, where $g_N(v) \in \mathcal{T}_N$ is a fixed descendant of $v \in 
\Tilde{\mathcal{T}}^{~}_{\lfloor \frac{N}{2} \rfloor}.$ 
We also set $\tilde{Z}_N := |\tilde{\mathcal{T}}_N|.$ 
By Lemma \ref{scgwlem} (1), 
$(\tilde{Z}_N)_{N \ge 0}$ is a Galton-Watson process with mean $m$ and zero extinction probability.
By applying Lemma \ref{scgwlem} (2) to $(\tilde{Z}_N)_{N \ge 0},$ we have
$$\lim_{N \to \infty} \frac{\tilde{Z}_{N+1}}{\tilde{Z}_N} = m, \textbf{P} \text{-a.s., and so} \lim_{N \to \infty} (\tilde{Z}_N)^{1/N} = m, \textbf{P} \text{-a.s.}$$
In particular, we have 
$|V^{'}| = \tilde{Z}_{\lfloor \frac{N}{2} \rfloor}  \ge {\alpha}^{\lfloor \frac{N}{2} \rfloor}$ for sufficiently large $N \in \mathbb{N},$ $\textbf{P}$-a.s., for all $1 < \alpha < m.$
We also know that $R_{\text{eff} }^N(x, y) > 2 \lfloor \frac{N}{2} \rfloor$ for all $x, y \in V^{'}, x \neq y,$ $\textbf{P}$-a.s.
Therefore, $\{B_{\text{eff}}^N (x, \lfloor \frac{N}{2} \rfloor) : x \in V^{'} \}$ is a packing for $\mathcal{T}_{\le N}$ 
and $\log \{n_{\text{pac}} (\mathcal{T}_{\le N}, \lfloor \frac{N}{2} \rfloor) \} \ge \lfloor \frac{N}{2} \rfloor \log \alpha,$ for sufficiently large $N \in \mathbb{N},$ $\textbf{P}$-a.s.,
for all $1 < \alpha < m.$ 
By Remark \ref{type2rems} (2), the conclusion holds. $\Box$ \\

\subsection{The largest supercritical percolation cluster inside a box in $\mathbb{Z}^d$}
\label{sec:8}
~~~~We consider Bernoulli bond percolation model on $\mathbb{Z}^d$. 
In this model, each edge in $\mathbb{E}^d$ is open with probability $p$ and closed with probability $1 - p$ independently, where
$\mathbb{E}^d := \{ \{x, y \} : x, y \in \mathbb{Z}^d, \sum_{i = 1}^d |x_i - y_i| = 1 \}$ and $x_i$ is the $i$ th coordinate of $x \in \mathbb{Z}^d$.
We write the corresponding probability measure on $\{0, 1 \}^{\mathbb{E}^d}$ by $\textbf{P}_p$.
A sequence $\Gamma = (x^0, \dotsc, x^n)$ is an open path in $S \subset \mathbb{Z}^d$ connecting $x$ and $y$ if $x^0 = x, x^n = y, x^i \in S$ for all $0 \le i \le n$
and $\{x^{i - 1}, x^i \}$ is an open edge for all $1 \le i \le n$.
We define the cluster at $x$ in $S \subset \mathbb{Z}^d$ by
$$\mathcal{C}^S (x) := \{y \in S :~\text{there exists an open path in $S$ connecting $x$ and $y$} \}.$$  
The critical probability is defined by
$$p_c(\mathbb{Z}^d ) := \inf \{p : ~\textbf{P}_p (\mathcal{C}^{\mathbb{Z}^d} (0) ~\text{is infinite} ) > 0 \}.$$
Let $\mathcal{C}_d(N)$ be the largest cluster in a box $[-N, N]^d$.  We prove the following results.
\begin{prop} \label{scpcprop}
(1) For $d = 2, p > p_c(\mathbb{Z}^2)$, there exist $c_1, c_2 > 0$ such that 
$$\displaystyle \lim_{N \to \infty} \textbf{P}_p (c_1 N^2 (\log N)^2 \le t_{cov}(\mathcal{C}_2(N)) \le c_2 N^2 (\log N)^3) = 1.$$ 
(2) For $d \ge 3, p > p_c(\mathbb{Z}^d)$, there exist $c_3, c_4 > 0$ such that
$$\displaystyle \lim_{N \to \infty} \textbf{P}_p (c_3 N^d \log N \le t_{cov}(\mathcal{C}_d(N)) \le c_4 N^d (\log N)^{\frac{2d - 1}{d - 1} }) = 1.$$ 
\end{prop}
\begin{rem} \label{remscpc}
Unfortunately, we are not able to obtain the correct order of the cover time.
If $\text{diam}_R (\mathcal{C}_2 (N))$ is of order $\log N$ as stated in Corollary 3.1 of \cite{BK}, we can obtain the correct order $(N^2 (\log N)^2)$
of the cover time for $\mathcal{C}_2 (N).$ However, from the proof of Corollary 3.1 of \cite{BK}, we can only obtain that $\text{diam}_R (\mathcal{C}_2 (N))$ is of order
$(\log N)^2.$ In particular, we can only state that $t_{\text{cov}} (\mathcal{C}_2 (N))$ is of order $N^2 (\log N)^3.$
\end{rem}
We use the following lemmas.
\begin{lem}(\cite{BM}, Proposition 1.2) \label{scpclem0}
For $d \ge 2, p > p_c(\mathbb{Z}^d),$ there exists $c > 0$ such that w.h.p.,
$$|\mathcal{C}_d(N)| \ge c N^d.$$ 
\end{lem}
Let $G = (V(G), E(G))$ be a finite graph. 
For $S \subset V(G),$ we define the external boundary of $S$ under the graph $G$ 
by $\partial_e S := \{x \in V(G) \backslash S : \text{there exists}~ y \in S ~\text{such that}~\{x, y \} \in E(G) \}.$
Set $L_x := \displaystyle \sum_{k = 1}^{\lfloor \log_2 |V(G)| \rfloor} \displaystyle \max \left(\frac{|S|}{|\partial_e S|^2} + \frac{1}{|\partial_e S|} \right),$ 
where the maximum is taken over all connected subsets $S$ of $V(G)$ satisfying
$x \in S$ and $|V(G)| / 2^{k + 1} < |S| \le |V(G)| / 2^k.$
\begin{lem}(\cite{BK}, Theorem 2.1) \label{scpclem1}
Let $G = (V(G), E(G))$ be a finite graph. There exists $c > 0$ such that for all $x, y \in V(G)$,
$$R_{\text{eff}} (x, y) \le c (L_x + L_y).$$
\end{lem}
\begin{lem}(\cite{Pe}, Corollary 1.4) \label{scpclem2}
Fix $d \ge 2, p > p_c (\mathbb{Z}^d).$ There exist $c, c^{'} > 0$ such that 
\begin{align*}
\lim_{N \to \infty} \textbf{P}_p \Big(|\partial_e S| \ge c |S|^{1 - 1/d} ~ &\text{for all connected subsets} ~S \subset \mathcal{C}_d(N)  \\
&\text{with} ~c^{'} (\log N)^{\frac{d}{d - 1}} \le |S| \le \frac{|\mathcal{C}_d(N)|}{2} \Big) = 1,
\end{align*}
where $\partial_e S$ is the external boundary of $S$ under the graph $\mathcal{C}_d(N).$
\end{lem}
\textit{Proof of Proposition \ref{scpcprop}.}
First, we prove the upper bounds by checking (\ref{condition-of-type1-0}) in Theorem \ref{type1} with $\log v(N) = \log N$ and $r(N) = (\log N)^{\frac{d}{d - 1}}$. 
It is clear that $|\mathcal{C}_d(N)| \le |[-N, N]^d \cap \mathbb{Z}^d| \le (2N + 1)^d,$ $\textbf{P}$-a.s. 
If $|\partial_e S| \ge c |S|^{1 - 1/d}$ for all connected subset $S \subset \mathcal{C}_d(N)$ with 
$c^{'} (\log N)^{\frac{d}{d - 1}} \le |S| \le \frac{|\mathcal{C}_d(N)|}{2},$
then we get for some $c_5 > 0,$
\begin{align*}
&\sum_{k = 1}^{\lfloor \log_2 |\mathcal{C}_d(N)| \rfloor} \max \Big \{ \frac{|S|}{|\partial_e S|^2} + \frac{1}{|\partial_e S|} 
: S~\text{is a connected subset of}~\mathcal{C}_d(N) \\
&~~~~~~~~~~~~~~~~~~~~~\text{satisfying}~x \in S ~\text{and}~ |\mathcal{C}_d(N)| / 2^{k + 1} < |S| \le |\mathcal{C}_d(N)| / 2^k \Big \} \\
&\le \sum_{k = 1}^{\lfloor \log_2 \{|\mathcal{C}_d(N)| / c^{'} (\log N)^{\frac{d}{d - 1}} \} \rfloor - 1} \Big (\frac{1}{c^2} + \frac{1}{c} \Big) \\
&+ \sum_{k = \lfloor \log_2 \{|\mathcal{C}_d(N)| / c^{'} (\log N)^{\frac{d}{d - 1}} \} \rfloor}^{\lfloor \log_2 |\mathcal{C}_d(N)| \rfloor} \Big (\frac{|\mathcal{C}_d (N)|}{2^k} + 1 \Big ) \\
&\le c_5 (\log N)^{\frac{d}{d - 1}} ~\text{for all}~x \in \mathcal{C}_d (N). 
\end{align*}
Therefore, by Lemma \ref{scpclem1} and Lemma \ref{scpclem2}, there exists $c_6 > 0$ such that w.h.p.,
$$\text{diam}_R(\mathcal{C}_d(N)) \le c_6 (\log N)^{\frac{d}{d - 1}}.$$
By Theorem \ref{type1} (1), we obtain the upper bound.

Next, we prove the lower bound for $d = 2$ by checking (\ref{condition-of-type1}) in Theorem \ref{type1}
with $\log v(N) = \log N$ and $r(N) = \log N.$ \\
If $|\mathcal{C}_2(N)| \ge c_7 N^2$, there exist $c_8 > 0, x, y \in \mathcal{C}_2(N)$ such that $d_{\mathbb{Z}^2}(x,y) > c_8 N.$ \\
We define a square with side length $2k$ centered at $u$ and its internal boundary by 
$$Q(u, k) := \{v \in \mathbb{Z}^2 : v_i \in [u_i - k, u_i + k], i = 1, 2 \},$$
$$\partial_i Q(u, k ) := \{v \in Q(u, k) : \exists w \in \mathbb{Z}^2 \backslash Q(u, k) ~\text{such that}~ \{v, w \} \in \mathbb{E}^2 \}.$$
Since $y \notin Q(x, \lfloor \frac{c_8}{2} N \rfloor ),$ there exists $x^k \in \mathcal{C}_2(N)$ such that
$x^k \in \partial_i Q(x, k\lfloor \sqrt{N} \rfloor )$ for all $0 \le k \le \frac{\lfloor \frac{c_8}{2} N \rfloor }{\lfloor \sqrt{N} \rfloor}.$ 
Fix $x^k, x^{\ell}, 0 \le k < \ell \le \frac{\lfloor \frac{c_8}{2} N \rfloor }{\lfloor \sqrt{N} \rfloor}.$
Since $d_{\mathbb{Z}^2 }(x^k, x^{\ell} ) \ge \lfloor \sqrt{N} \rfloor$, there exists a positive integer $a(N) \in [\lfloor \frac{\lfloor \sqrt{N} \rfloor}{2} \rfloor, \infty)$ such that
$x^{\ell} \in \partial_i Q(x^k, a(N) ).$
We write $\Pi_j := \{\{u, v \} \in \mathbb{E}^2 : u \in \partial_i Q(x^k, j - 1) ~\text{and}~v \in \partial_i Q(x^k, j) \}, 1 \le j \le a(N).$
Under the induced graph $G_{cN}$ with vertex set $[- cN, cN]^2 \cap \mathbb{Z}^2$ for some sufficiently large constant $c > 0,$
$(\Pi_j)_{1 \le j \le a(N)}$ is a sequence of edge-cutsets between $x^k$ and $x^{\ell}$. 
So, we have by Lemma \ref{Nash-Williams} that for some $ c_9 > 0$,
\begin{equation} \label{scpcineq}
R_{\text{eff}}^N (x^k, x^{\ell}) \ge R_{\text{eff}}^{G_{cN}} (x^k, x^{\ell}) \ge c_9 \log N,
\end{equation}
where $R_{\text{eff}}^{G_{cN}} (\cdot, \cdot)$ is the effective resistance in the graph $G_{cN}.$ \\
Set $V^{'} := \{ x^0, x^1, \dotsc, x^{\lfloor \frac{\lfloor \frac{c_8}{2} N \rfloor }{\lfloor \sqrt{N} \rfloor} \rfloor} \}.$
By (\ref{scpcineq}), $\{B_{\text{eff}}^N (x, \frac{c_9}{4} \log N) : x \in V^{'} \}$ is a packing for $\mathcal{C}_2(N).$
So, there exists $c_{10} > 0$ such that w.h.p.,
$$\log \{n_{\text{pac}} (\mathcal{C}_2(N), \frac{c_9}{4} \log N) \} \ge c_{10} \log N.$$ 
Therefore, by Theorem \ref{type1} (2) and Lemma \ref{scpclem0}, we get the lower bound for $d = 2.$

We next prove the lower bound for $d \ge 3$ by checking (\ref{condition-of-type1}) in Theorem \ref{type1}
with $\log v(N) = \log N$ and $r(N) = 1.$
Fix $u, v \in \mathcal{C}_d (N), u \neq v.$ Set $\Pi := \{ \{u, x \} : \{u, x \} \in E(\mathcal{C}_d(N)) \}.$ $\Pi$ is an edge-cutset between $u$ and $v$
in the graph $\mathcal{C}_d (N).$ So, by Lemma \ref{Nash-Williams}, we have that $R_{\text{eff}}^N (u, v) \ge 1/|\Pi| \ge 1/2d.$ 
In particular, $\{B_{\text{eff}}^N (x, 1/8d) : x \in \mathcal{C}_d(N) \}$ is a packing for $\mathcal{C}_d(N).$
So, by Lemma \ref{scpclem0}, we have for some $c_{11} > 0,$
$$\log \{n_{\text{pac}} (\mathcal{C}_d (N), 1/8d) \} \ge c_{11} \log N ~\text{w.h.p.}$$
Therefore, by Theorem \ref{type1} (2) and Lemma \ref{scpclem0}, we obtain the lower bound for $d \ge 3$. $\Box$ \\

\subsection{Supercritical Erd\H{o}s-R\'{e}nyi random graph I}
\label{sec:9}
~~~~Let $G(N, p)$ be the Erd\H{o}s-R\'{e}nyi random graph. This is obtained from the complete graph with $N$ vertices by retaining each edge with probability $p$ independently. We 
assume that $p = \frac{c}{N}$, where $c > 1$ is a positive constant. Let $\mathcal{C}^N$ be the largest connected component of $G(N, p)$.

We revisit Theorem 2a of \cite {CF1}. Note that Cooper and Frieze \cite{CF1} has obtained a better estimate than the following Proposition \ref{scerprop}. 
See Remark \ref{scerrem} below.
\begin{prop} \label{scerprop} There exist $c_1, c_2 > 0$ such that
$$\displaystyle \lim_{N \to \infty} \textbf{P} (c_1 N (\log N)^2 \le t_{cov}(\mathcal{C}^N) \le c_2 N (\log N)^2) = 1,$$ 
and $(t_{cov}(\mathcal{C}^N))_{N \in \mathbb{N}}$ is Type 1.
\end{prop}
\textit{Proof.}
We check (\ref{condition-of-type1-0}) and (\ref{condition-of-type1}) in Theorem \ref{type1} with $v(N) = N$ and $r(N) = \log N.$
It is known that w.h.p., $(1 - \epsilon) \frac{cx(2-x)}{2} N \le |E(\mathcal{C}^N)| \le (1 + \epsilon) \frac{cx(2-x)}{2} N$ for any $\epsilon > 0,$
where $x$ is the solution of $x = 1 - e^{-cx} $ in $(0, 1)$ (see Section 3.1.3 of \cite{CF1}).
By Theorem 6 of \cite{CL}, there exists $c_3 > 0$ such that w.h.p.,
$$\text{diam}_R(\mathcal{C}^N) \le \text{diam}(\mathcal{C}^N) \le c_3 \log N.$$
The largest connected component $\mathcal{C}^N$ consists of a 2-core $C_2$ (the largest subgraph of $\mathcal{C}^N$ with minimum degree $2$) and
a mantle $\textbf{M}$ (a collection of trees which are sprouting from different vertices of $C_2$). 
By Lemma 9 and \textbf{P7a} of \cite{CF1}, w.h.p., there exists a subset $V^{'} \subset \mathcal{C}^N$ which satisfies the following:\\ \\
(i) Every $v \in V^{'}$ is a leaf of a tree $T_v$ in $\textbf{M},$\\
(ii) Let $w(v)$ be the root of $T_v$. Then, $d_{\mathcal{C}^N} (v, w(v)) = \lceil \log N/(2(cx - \log c)) \rceil$,\\
(iii) For any $\epsilon > 0, c_4 N^{1/2 - \epsilon} \le |V^{'}| \le c_5 N^{1/2 + \epsilon}, ~T_u \neq T_v ~\text{for all}~ u, v \in V^{'}, u \neq v,$ 
where $c_4, c_5 > 0$ are some constants. \\ \\
(Indeed, choose `special vertices' in their terminology in Section 3.1.2 of \cite{CF1}.) \\
In particular, if $u, v \in V^{'}, u \neq v$, then every path from $u$ to $v$ contains a common path of length $\lceil \log N/(2(cx - \log c)) \rceil$.
By Lemma \ref{Nash-Williams}, we have that $R_{\text{eff} }^N(u, v) >  2 \lceil \log N/(2(cx - \log c)) \rceil.$
In particular, $\{B_{\text{eff}}^N (v, \lceil \frac{\log N}{2(cx - \log c)} \rceil) : v \in V^{'} \}$ is a packing for $\mathcal{C}^N.$
So, we have for some $c_6 > 0,$
$$\log \Big \{n_{\text{pac}} \Big(\mathcal{C}^N, \Big \lceil \frac{\log N}{2(cx - \log c)} \Big \rceil \Big) \Big \} \ge \log |V^{'}| \ge c_6 \log N,~\text{w.h.p.}~\Box$$

\begin{rem} \label{scerrem}
In \cite{CF1}, Cooper and Frieze proved that for any $\epsilon > 0,$ w.h.p.,
$$(1 - \epsilon) \frac{cx(2 - x)}{4(cx - \log c)} N(\log N)^2 \le t_{\text{cov}}(\mathcal{C}^N) \le (1 + \epsilon) \frac{cx(2 - x)}{4(cx - \log c)} N(\log N)^2.$$ 
\end{rem}
\subsection{Supercritical Erd\H{o}s-R\'{e}nyi random graph II}
\label{sec:10}
~~~~We consider the Erd\H{o}s-R\'{e}nyi random graph $G(N, p)$ again. 
Here we assume that $p = \frac{f(N)}{N}$, where $\lim_{N \to \infty} \frac{\log N}{f(N)} = \lim_{N \to \infty} \frac{f(N)}{N^{1/2}} = 0.$
In this regime, $G(N, p)$ is connected w.h.p.

We revisit Theorem 1.1(i) of \cite{J}. Note that Jonasson \cite{J} has obtained a better estimate than the following Proposition \ref{scer2prop}. 
See Remark \ref{scer2rem} below.
\begin{prop} \label{scer2prop}
There exist $c_1, c_2 > 0$ such that
$$\displaystyle \lim_{N \to \infty} \textbf{P} (c_1 N \log N \le t_{cov}(G(N, p)) \le c_2 N \log N ) = 1,$$
and $(t_{cov}(G(N, p)))_{N \in \mathbb{N}}$ is Type 1.
\end{prop}
In the proof, we use the following lemma. 
\begin{lem} (Lemma 3.1 and Proposition 3.1 of \cite{J}) \label{scer2lem}
Fix any $\epsilon > 0$. Then, w.h.p.,
$$(1 - \epsilon) f(N) \le \mu_x^N \le (1 + \epsilon) f(N), ~\text{for all}~x \in G(N, p), $$
$$\text{diam}_R (G(N, p)) \le \frac{2}{(1 - \epsilon) f(N)}.$$
\end{lem}
\textit{Proof of Proposition \ref{scer2prop}.}
We check (\ref{condition-of-type1-0}) and (\ref{condition-of-type1}) in Theorem \ref{type1} with $v(N) = f(N) N$ and $r(N) = 1/f(N).$
By Lemma \ref{scer2lem}, there exist $c_3, c_4 > 0$ such that w.h.p., 
$$c_3 f(N) N \le |E(G(N, p))| \le c_4 f(N) N.$$
By this together with Lemma \ref{scer2lem}, (\ref{condition-of-type1-0}) holds. \\
By Lemma \ref{Nash-Williams} and Lemma \ref{scer2lem}, there exists $c_5 > 0$ such that $\{B_{\text{eff}}^N (x, \frac{c_5}{f(N)}) : x \in G(N, p) \}$ is a packing for $G(N, p)$ w.h.p.
So, we have for some $c_6 > 0,$
$$\log \{n_{\text{pac}} (G(N, p), c_5 / f(N) ) \} \ge c_6 \log \{ f(N)N \},~\text{w.h.p.}~\Box$$

\begin{rem} \label{scer2rem}
In \cite{J}, Jonasson proved that if $\lim_{N \to \infty} \frac{\log N}{f(N)} = 0$, then for any $\epsilon > 0,$ w.h.p., 
$$(1 - \epsilon) N \log N \le t_{\text{cov}}(G(N, p)) \le (1 + \epsilon) N \log N.$$
\end{rem}

\subsection{The incipient infinite cluster for critical Galton-Watson family trees}
\label{sec:11}
~~~~Let $(Z_N)_{N \ge 0}$ be a critical Galton-Watson process with offspring distribution $Z$ in the domain of attraction of a stable law with index $\alpha \in (1, 2 ].$
That is, there exists a sequence $(a_N)_{N \ge 0}$ such that $\frac{Z[N] - N}{a_N} \stackrel{d}{\to} X$,
where $\text{E}e^{- \lambda X} = e^{- \lambda^{\alpha}}$ and $Z[N]$ is the sum of $N$ i.i.d copies of $Z$. 
We write $\mathcal{T}$ to denote its family tree.
We use the notation $\mathcal{T}_{\le N}, \mathcal{T}_N$ as in Subsection 3.1. 
We set $p_N := \text{P} (Z_N > 0)$. 
In \cite{Ke}, Kesten considered the Galton-Watson tree conditioned to survive:

\begin{lem} (\cite{Ke}, Lemma 1.14) For any family tree $T$ of $k$ generations,
$$\lim_{N \to \infty} \text{P} (\mathcal{T}_{\le k} = T | Z_N > 0) = |T_k| \text{P} (\mathcal{T}_{\le k} = T ).$$
We set $P_0 (T) = |T_k| P(\mathcal{T}_{\le k} = T)$. $P_0$ has a unique extension to a probability measure $\textbf{P}$ on the set of infinite family trees. \\
\end{lem} 
By this lemma, we can take a family tree with the distribution $\textbf{P}$. We write this by $\mathcal{T}^{*}$ and call it incipient infinite cluster. 
We set $Z_N^{*} := |\mathcal{T}_N^{*}|.$
\begin{prop} \label{critical GW} There exist $c_1, c_2, c > 0$ such that for all $\lambda, N \ge c$,
$$\textbf{P}(t_{\text{cov}}(\mathcal{T}_{\le N}^{*}) \ge \lambda N^{\frac{2 \alpha - 1}{\alpha - 1}} \ell(N)^{- 1} ) \le c_1 \lambda^{- c_2}, $$
$$\textbf{P}(t_{\text{cov}}(\mathcal{T}_{\le N}^{*}) \le \lambda^{- 1} N^{\frac{2 \alpha - 1}{\alpha - 1}} \ell(N)^{- 1}) \le c_1 \lambda^{- c_2},$$
where $\ell (N)$ is a slowly varying function at infinity satisfying $p_N = N^{- \frac{1}{\alpha - 1}} \ell (N).$ \\
Furthermore, $(t_{\text{cov}}(\mathcal{T}_{\le N}^{*}))_{N \in \mathbb{N}}$ is Type 2.
\end{prop}
\begin{rem}
Barlow, Ding, Nachmias and Peres \cite{BDNP} proved that in the case $\alpha = 2$, conditioned on the event $\{|\mathcal{T}| \in [N, 2N] \},$
$t_{\text{cov}} (\mathcal{T}) / N^{\frac{3}{2}}$ is tight.
\end{rem}
In the proof, we use the following facts.
\begin{lem} (Proposition 2.2, 2.5, 2.7 and Lemma 2.3 of \cite{CK}) \label{cgwlem} \\
(1) There exists a slowly varying function at infinity $\ell(N)$ which satisfies that $p_N = N^{- \frac{1}{\alpha - 1}} \ell(N)$ 
and that for any $\epsilon > 0$, there exist $c_3, c_4 > 0$ such that
\begin{equation*}
c_3 \left(\frac{N}{N^{'}} \right)^{- \epsilon} \le \frac{\ell(N)}{\ell(N^{'})} \le c_4 \left(\frac{N}{N^{'}} \right)^{\epsilon}, ~\text{for all}~1 \le N^{'} \le N.
\end{equation*} 
(2) Set $J(\lambda) := \{N \in \mathbb{N}: Z_N^{*} \le \lambda p_N^{- 1}, |E(\mathcal{T}_{\le N}^{*} )| \ge
\lambda^{- 1} N p_N^{- 1}, |\mathcal{T}_{\le N}^{*} | \le \lambda N p_N^{- 1} \}.$ 
Then there exist $c_5, c_6 > 0$ such that for all $N \in \mathbb{N}$ and $\lambda > 0$,
$$\textbf{P} (N \in J(\lambda) ) \ge 1 - c_5 \lambda^{- c_6}.$$
\end{lem}
\textit{Proof of Proposition \ref{critical GW}.} \\
By Lemma \ref{cgwlem} (2) and the fact that $N \le \text{diam}_R (\mathcal{T}_{\le N}^{*}) \le 2N$ $\textbf{P}$-a.s., 
the conditions (\ref{condition-of-type2-(1)-1}) and (\ref{condition-of-type2-(3)}) in Theorem \ref{type2} hold for $v(N) = Np_N^{-1}$ and $r(N) = N.$ 
So, we only need to check (\ref{condition-of-type2-(1)-2}) with $r(N) = N.$ \\
The idea of the following argument came from the proof of Theorem 3.1 of \cite{BDNP}.
We write $\mathcal{T}^{*, x}$ to denote the subtree rooted at $x \in \mathcal{T}^{*}$.
Set $r_{k,j}^N := \lfloor \frac{j}{2^{k+2}} N \rfloor, k \in \mathbb{N}, 0 \le j \le 2^{k+2}.$  \\
Fix $k \in \mathbb{N}$ and $0 \le j \le 2^{k+2} - 1.$
We say that $x \in \mathcal{T}_{r_{k,j}^N }^{*} $ is $k$-good if $\mathcal{T}_{(r_{k, j + 1}^N - r_{k,j}^N  ) }^{*, x} \neq \emptyset.$
We assume $\lambda \ge c_7$, where $c_7$ is a sufficiently large positive constant. 
Set for all $0 \le j \le 2^{k + 2} - 1,$
$$A_{k, j}^N := \{x \in \mathcal{T}_{r_{k,j}^N}^{*} : x~\text{is $k$-good} \}.$$
We define

$$A_{k}^N := \begin{cases}
                           \displaystyle \bigcup_{j = 0}^{2^{k+2} - 1} A_{k, j}^N & \text{if} ~0 \le k \le \lfloor \frac{\log N}{\log 2} \rfloor - 2, \\ 
                             \\
                            \mathcal{T}_{\le N}^{*} & \text{otherwise}.
                            \end{cases}$$
We define $\ell_k^N := \frac{\text{diam}_R (\mathcal{T}_{\le N}^{*})}{2^k}$ for $0 \le k \le \lfloor \frac{\log N}{\log 2} \rfloor - 2$ and $\ell_k^N = 0$ otherwise. \\
Since $\{B_{\text{eff}}^N (x, \ell_k^N) : x \in A_k^N \}$ is a covering for $\mathcal{T}_{\le N}^{*}$ for all $k \ge 0,$ we get for all $k \ge 0,$
\begin{equation}
n_{\text{cov}} (\mathcal{T}_{\le N}^{*}, \ell_k^N) \le |A_k^N|. \label{covering-number-ct}
\end{equation}
Fix $0 \le k \le \lfloor \frac{\log N}{\log 2} \rfloor - 2$ and $1 \le j \le 2^{k+2} - 1.$
By Lemma 2.2 of \cite{Ke} (note that in \cite{Ke}, Kesten assumed the variance of offspring distribution is finite, but the same result holds under our situation), 
for $\tilde{\lambda} > 0$, 
\begin{align*}
&~~~~\textbf{P}(|A_{k, j}^N| \ge \tilde{\lambda} | \mathcal{T}_{\le r_{k,j}^N}^{*} = T, H_{\le r_{k,j}^N} = (v_i)_{0 \le i \le r_{k,j}^N} ) \\
&= \textbf{P}(|A_{k, j}^N \backslash \{ v_{r_{k,j}^N} \}| \ge \tilde{\lambda} - 1 | \mathcal{T}_{\le r_{k,j}^N}^{*} = T, H_{\le r_{k,j}^N} = (v_i)_{0 \le i \le r_{k,j}^N} ) \\
&= \textbf{P}(\text{Bin}(|T_{r_{k,j}^N}| - 1, p_{(r_{k,j + 1}^N - r_{k,j}^N)} ) \ge \tilde{\lambda} - 1 ), 
\end{align*}
where $T$ is a family tree of $r_{k,j}^N$ generations, $H_{\le r_{k,j}^N}$ is a backbone (the unique infinite line of descent of $\mathcal{T}^{*}$)
up to $r_{k,j}^N$ th level and $(v_i)_{0 \le i \le r_{k,j}^N}$ is a sequence of vertices
such that $v_i \in T_i$ for all $0 \le i \le r_{k,j}^N.$
We also note that for all $0 \le m \le \Big \lfloor \frac{\tilde{\lambda}}{2p_{(r_{k,j + 1}^N - r_{k,j}^N )}} \Big \rfloor,$
\begin{align*}
&\textbf{P}(\text{Bin}(m, p_{(r_{k,j + 1}^N - r_{k,j}^N)} ) \ge \tilde{\lambda} - 1) \\
&\le \textbf{P}(\text{Bin}(\Big \lfloor \frac{\tilde{\lambda}}{2p_{(r_{k,j + 1}^N - r_{k,j}^N )} } \Big \rfloor, p_{(r_{k,j + 1}^N - r_{k,j}^N)} ) \ge \tilde{\lambda} - 1 ).
\end{align*}
Therefore, for $\tilde{\lambda} > 2$,
\begin{align*}
&~~~~\textbf{P}(|A_{k, j}^N| \ge \tilde{\lambda} ) \\
&\le \textbf{P}(\text{Bin}(\Big \lfloor \frac{\tilde{\lambda}}{2p_{(r_{k,j + 1}^N - r_{k,j}^N)}} \Big \rfloor, p_{(r_{k,j + 1}^N - r_{k,j}^N)} ) \ge \tilde{\lambda} - 1 ) \\
&+ \textbf{P}(Z_{r_{k,j}^N}^{*} > \Big \lfloor \frac{\tilde{\lambda}}{2p_{(r_{k,j + 1}^N - r_{k,j}^N)}} \Big \rfloor ).
\end{align*}         

By the Chebyshev inequality, the first term is bounded by $\frac{2 \tilde{\lambda}}{(\tilde{\lambda} - 2 )^2}$.
By Lemma \ref{cgwlem} (1) (2), the second term is bounded by $c_8 j^{c_9} \tilde{\lambda}^{- c_{10}}$
for some $c_8, c_9, c_{10} > 0.$
So, we have that
\begin{align*}
&~~~~\textbf{P}(|A_k^N| \ge \exp (\lambda 2^{k/2}) ) \\
&\le \textbf{P} \bigg(\displaystyle \bigcup_{j = 1}^{2^{k+2} - 1} \bigg \{ |A_{k, j}^N| \ge \frac{\exp (\lambda 2^{k/2}) - 1}{2^{k+2}} \bigg \} \bigg) \\                   
&\le \displaystyle \sum_{j = 1}^{2^{k+2} - 1} 
\Bigg \{\frac{2 \cdot \frac{\exp (\lambda 2^{k/2}) - 1}{2^{k+2}}}{(\frac{\exp (\lambda 2^{k/2}) - 1}{2^{k+2}} - 2)^2 }
+ c_8 j^{c_9} \bigg(\frac{\exp (\lambda 2^{k/2}) - 1}{2^{k+2}} \bigg)^{- c_{10} } \Bigg \} \\                 
&\le c_{11} 2^{-k} \lambda^{- c_{12}}  ~~~~\text{for some}~ c_{11}, c_{12} >0.
\end{align*}
From this fact, we have that 
\begin{equation}
\textbf{P} \Bigg(\displaystyle \bigcup_{k = 0}^{\lfloor \frac{\log N}{\log 2} \rfloor - 2} \bigg \{ |A_k^N| \ge \exp (\lambda 2^{k/2}) \bigg \} \Bigg) 
\le 2c_{11} \lambda^{- c_{12}}. \label{covering-number-ct-2}
\end{equation}
If $|A_k^N| \le \exp (\lambda 2^{k/2})$ for all $0 \le k \le \lfloor \frac{\log N}{\log 2} \rfloor - 2$ and $|\mathcal{T}_{\le N}^{*}| \le \lambda N p_N^{-1},$ we have 
by (\ref{covering-number-ct}),
$$\sum_{k = 1}^{\lfloor \frac{\log N}{\log 2} \rfloor - 1} \sqrt{\ell_{k-1}^N \log \{n_{\text{cov}} (\mathcal{T}_{\le N}^{*}, \ell_k^N) \}} \le c_{13} \sqrt{\lambda N}$$ 
for some $c_{13} > 0$. \\
So, by (\ref{covering-number-ct-2}) and Lemma \ref{cgwlem} (2), (\ref{condition-of-type2-(1)-2}) in Theorem \ref{type2} holds with $r(N) = N.$ $\Box$ \\

We can also say that $t_{\text{cov}}(\mathcal{T}_{\le N}^{*}) N^{- \frac{2 \alpha - 1}{\alpha - 1}} \ell (N)$ is not concentrated. 
\begin{prop} \label{cgwprop2}
For all $\lambda \ge 1,$
$$\liminf_{N \to \infty} \textbf{P}(t_{\text{cov}}(\mathcal{T}_{\le N}^{*}) N^{- \frac{2 \alpha - 1}{\alpha - 1}} \ell (N) \ge \lambda ) > 0.$$
\end{prop}
To prove this fact, we use the following result.
\begin{lem} (\cite{Pa}, Theorem 4) \label{pakes}
The random variable $Z_N^{*} p_N$ converges in law to a random variable $Z^{*}$ 
with $\mathbb{E} (e^{- \theta Z^{*}}) = (1 + \theta^{\alpha - 1})^{- \frac{\alpha}{\alpha - 1}}$ for $\theta \ge 0.$
\end{lem}
\textit{Proof of Proposition \ref{cgwprop2}.} \\
By the fact that $t_{\text{cov}} (\mathcal{T}_{\le N}^{*}) \ge t_{\text{hit}} (\mathcal{T}_{\le N}^{*}) \ge \frac{1}{2} N |E(\mathcal{T}_{\le N}^{*})|$ (we have 
used Lemma \ref{commute time}), 
for $\lambda >0,$
$$\textbf{P}(t_{\text{cov}}(\mathcal{T}_{\le N}^{*}) N^{- \frac{2 \alpha - 1}{\alpha - 1}} \ell (N) \ge \lambda ) \ge \textbf{P}(|E(\mathcal{T}_{\le N}^{*})| \ge 2 \lambda N p_{N}^{- 1} ).$$
Using the proof of Proposition 2.5 of \cite{CK} (in page 1429) when $\alpha \in (1, 2)$ and Lemma \ref{pakes} when $\alpha = 2$,
we have that for $\lambda \ge 1$ and some $c_{14}, c_{15} > 0,$
\begin{equation*}
\liminf_{N \to \infty} \textbf{P} (|E(\mathcal{T}_{\le N}^{*})| \ge \lambda N p_N^{- 1}) \ge c_{14} \liminf_{N \to \infty} \textbf{P} (Z_{N^{'}}^{*} p_{N^{'}} > c_{15} \lambda) > 0,
\end{equation*}
where $N^{'} = \lfloor \frac{N}{3} \rfloor.$
This implies the conclusion. $\Box$

\subsection{Critical percolation clusters}
\label{sec:12}
~~~~Let $\Hat{G}^N$ be a graph with $N$ vertices and the maximum degree $d \in [3, N - 1]$. $\Hat{G}_{p}^N$ is obtained by retaining each edge of $\Hat{G}^N$ with probability $p \in (0, 1)$ independently. Let  $\mathcal{C}^N$ be the largest connected component of $\Hat{G}_p^N$ and $\mathcal{C}(x)$ be the connected component of $\Hat{G}_p^N$ which contains $x \in V(\Hat{G}^N)$. We define balls and their boundaries as follows:
$$B_p(x, r; \Hat{G}^N) := \{ y \in V(\Hat{G}_p^N) : d_{\Hat{G}_p^N} (x, y) \le r \},$$
$$\partial B_p (x, r; \Hat{G}^N) := \{y \in V(\Hat{G}_p^N) : d_{\Hat{G}_p^N} (x, y) = r \}.$$
We also set
$$\Gamma_p(x, r; \Hat{G}^N ) := \sup_{G \subset \Hat{G}^N} \textbf{P}_{G} (H_p(x, r; G) ),$$
where $H_p (x, r; G ) := \{\partial B_p (x, r; G) \neq \emptyset \},$ 
the supremum is taken over all subgraphs of $\Hat{G}^N$ and $\textbf{P}_G$ is a percolation probability measure on $G$.
In particular, we write $\textbf{P} := \textbf{P}_{\Hat{G}^N}$. \\
We assume that
\begin{equation}
p \le \frac{1 + A N^{- 1/3} }{d - 1} ~\text{for some}~ A \in \mathbb{R}, \label{ineqcr1}
\end{equation}
and that there exist $c_1, c_2 > 0$ and $a : (0, \infty) \to (0, \infty)$
such that for sufficiently large $\lambda > 0$ and $N \ge a(\lambda),$
\begin{equation}
\textbf{P}(|\mathcal{C}^N| \le {\lambda}^{- 1} N^{\frac{2}{3} } ) \le c_1 {\lambda}^{- c_2}.  \label{ineqcr2}
\end{equation}
\begin{rem}
In the case that $\Hat{G}^N$ is the complete graph with $N$ vertices and $p = 1/N$, it is known that (\ref{ineqcr2}) holds (see Theorem 2 of \cite{NP2}).
\end{rem}
We revisit Theorem 3.1 of \cite{BDNP}. 
\begin{prop} \label{crprop}
Under the assumption (\ref{ineqcr1}) and (\ref{ineqcr2}),
there exist $c_3, c_4 > 0$ such that for sufficiently large $\lambda > 0$ and $N \ge \max \{\lambda^3, a(\lambda) \},$
$$\textbf{P}(t_\text{cov} (\mathcal{C}^N ) > \lambda N ) \le c_3 \lambda^{- c_4 }, ~\textbf{P}(t_\text{cov} (\mathcal{C}^N ) < {\lambda}^{- 1} N ) \le c_3 \lambda^{- c_4 },$$
and $(t_\text{cov} (\mathcal{C}^N) )_{N \in \mathbb{N}}$ is Type 2.
\end{prop}
\begin{rem} \label{crrem-0}
Barlow, Ding, Nachmias and Peres \cite{BDNP} have already considered the cover time for the critical random graphs.
\end{rem}
To prove this proposition, we use the following facts (most of them are proved in \cite{NP}).
\begin{lem} \label{crlem}
(1) There exists $c_5 > 0$ such that for all subgraphs $G \subset \Hat{G}^N, x \in V(G), \lambda >0$ and sufficiently large $N \in \mathbb{N},$ 
\begin{equation}
\textbf{E}|B_p(x, r; G )| \le c_5 e^{A \lambda} r,~~~~~~~ \text{for all}~ r \le \lambda N^{\frac{1}{3} }, \label{ineqcr3} 
\end{equation}
\begin{equation}
\Gamma_p (x, r; \Hat{G}^N) \le c_5/r, ~~~~~~~~~\text{for all} ~r \le N^{\frac{1}{3} }, \label{ineqprobab}
\end{equation}
where $A$ is the constant in (\ref{ineqcr1}). \\
(2) There exists $c_6 > 0$ such that for sufficiently large $N \in \mathbb{N}, \lambda > 0,$
$$\textbf{P} (|E (\mathcal{C}^N)| \ge \lambda N^{2/3} ) \le c_6 \lambda^{- 1},  \textbf{P} (\text{diam} (\mathcal{C}^N) \ge \lambda N^{1/3} ) \le c_6 \lambda^{- 1},$$
where $\text{diam} (\mathcal{C}^N) = \max_{x, y \in V(\mathcal{C}^N)} d_{\mathcal{C}^N} (x, y).$ \\
(3) There exists $c_7 > 0$ such that for sufficiently large $\lambda > 0$ and $N \ge \lambda^2,$
$$\textbf{P} (\exists x \in V(\Hat{G}^N), 
|\mathcal{C} (x) | > \lambda^{- 1/12} N^{2/3}~\text{and}~\text{diam}_R (\mathcal{C} (x)) < \lambda^{-1} N^{1/3} ) \le c_7 \lambda^{-1/6}.$$
(4) There exists $c_8 > 0$ such that for sufficiently large $\lambda > 0, N \ge \max \{\lambda^3, a(\lambda) \},$ 
$$\textbf{P} (|E (\mathcal{C}^N)| < \lambda^{-1} N^{2/3} ) \le c_8^{-1} \lambda^{- c_8},  \textbf{P} (\text{diam}_R (\mathcal{C}^N) < \lambda^{-1} N^{1/3} ) \le c_8^{-1} \lambda^{- c_8}.$$
\end{lem}
To prove (3) of this lemma, we use Proposition 5.6 in \cite{NP}. So, we recall some terms in \cite{NP}. \\
Fix $x \in \Hat{G}^N, r, L \in \mathbb{N}, k < r.$ 
For $j < r,$ a lane for $(x, r, j)$ is an edge $\{u, v \}$ with $u \in \partial B_p (x, j - 1; \Hat{G}_p^N)$ and $v \in \partial B_p (x, j ; \Hat{G}_p^N)$ 
such that a path from $u$ to a vertex in $\partial B_p (x, r ; \Hat{G}_p^N)$ passes $\{u,v \}$ and does not intersect $\partial B_p (x, j - 1; \Hat{G}_p^N)$ except the starting vertex. \\
We say $x$ is $L$-lane rich for $(k, r)$
if we have a subset $I \subset [\lfloor k/2 \rfloor, k] \cap \mathbb{Z}$ with $|I| > \lfloor \frac{1}{2} (k - \lfloor \frac{k}{2} \rfloor) \rfloor$
such that for any $j \in I$, there exist at least $L$ lanes for $(x, r, j).$ \\
\begin{lem} (\cite{NP}, Proposition 5.6) \label{l-lane-rich}
Suppose that $x \in V(\Hat{G}^N), L \in \mathbb{N}, k \le r/2$ and $r < N^{1/3}.$ Then there exists $c_9 > 0$ such that
$$\textbf{P} (x~\text{is $L$-lane rich for}~(k, r) ) \le c_9 L^{- 1} r^{- 1}.$$
\end{lem}
\textit{Proof of Lemma \ref{crlem}}
By the proof of Theorem 1.2 and Theorem 1.3 of \cite{NP} in page 1281, (1) holds.
The results of (2) are proved in \cite{NP} in page 1274 and 1283. \\
The result of (4) follows from (3) and  (\ref{ineqcr2}). 
So, we only prove (3). \\
We use Lemma \ref{l-lane-rich} with 
$$k = \Big \lfloor \frac{1}{2} \Big \lfloor N^{1/3} \Big (\frac{\lambda}{32} \Big)^{-1/3} \Big \rfloor \Big \rfloor,
 r = \Big \lfloor N^{1/3} \Big (\frac{\lambda}{32} \Big)^{-1/3} \Big \rfloor, L = \Big \lfloor \Big (\frac{\lambda}{32} \Big)^{2/3} \Big \rfloor.$$
Suppose that $x$ is not $L$-lane rich for $(k, r)$ and $\text{diam} (\mathcal{C} (x)) \ge (\frac{\lambda}{32})^{-1/3} N^{1/3}.$ \\
Since $x$ is not $L$-lane rich for $(k, r),$ there exists a subset $I \subset [\lfloor k/2 \rfloor, k] \cap \mathbb{Z}$ with $|I| \ge 
\lfloor \frac{1}{2} (k - \lfloor \frac{k}{2} \rfloor) \rfloor$ 
such that for all $j \in I,$ the number of lanes for $(x, r, j)$ is less than $L.$ \\
For $j \in I,$ let $\Pi_j$ be a set of all lanes for $(x, r, j).$ Note that by the property of $I,$ we have $|\Pi_j| \le L.$ \\
Because $\text{diam} (\mathcal{C} (x)) \ge (\frac{\lambda}{32})^{-1/3} N^{1/3},$ there exists a vertex $x_0$ in $\partial B_p (x, r; \Hat{G}_p^N).$
Since $\Pi_j$ is an edge-cutset between $x$ and $x_0$ for all $j \in I$, we get by Lemma \ref{Nash-Williams} for sufficiently large $\lambda > 0, N \ge \lambda^2,$
$$\text{diam}_R (\mathcal{C} (x)) \ge R_{\text{eff}}^N (x, x_0) \ge \sum_{j \in I} 1/|\Pi_j| \ge |I|/L \ge \lambda^{-1} N^{1/3}.$$
Therefore, we have
\begin{align*}
&\textbf{P} (|\mathcal{C} (x)| > \lambda^{-1/12} N^{2/3}, \text{diam}_R ( \mathcal{C} (x)) < \lambda^{-1} N^{1/3} ) \\
&\le \textbf{P} (|\mathcal{C} (x)| > \lambda^{-1/12} N^{2/3}, \text{diam} ( \mathcal{C} (x)) < \Big (\frac{\lambda}{32} \Big)^{-1} N^{1/3} ) \\
&+ \textbf{P} (x~\text{is $L$-lane rich for}~(k, r)). 
\end{align*}
By (\ref{ineqcr3}), we get for some $c_{10} > 0,$
\begin{align}
&\textbf{P} (|\mathcal{C} (x)| > \lambda^{-1/12} N^{2/3}, \text{diam} ( \mathcal{C} (x)) < \Big (\frac{\lambda}{32} \Big)^{-1} N^{1/3} ) \notag \\
&\le \textbf{P} (|B_p (x, r; \Hat{G}^N)| > \lambda^{-1/12} N^{2/3}) \notag \\
&\le \frac{\textbf{E} |B_p (x, r; \Hat{G}^N)|}{\lambda^{-1/12} N^{2/3}} \notag \\
&\le \frac{c_5 e^A r}{\lambda^{-1/12} N^{2/3}} \notag \\
&\le c_{10} \lambda^{- 1/4} N^{- 1/3}. \label{term1-2-cr}
\end{align}
By Lemma \ref{l-lane-rich} and (\ref{term1-2-cr}), we have for some $c_{11} > 0$ and sufficiently large $\lambda > 0, N \ge \lambda^2,$
\begin{equation}
\textbf{P} (|\mathcal{C} (x) | > \lambda^{- 1/12} N^{2/3}~\text{and}~\text{diam}_R (\mathcal{C} (x)) < \lambda^{-1} N^{1/3}) \le c_{11} \lambda^{- 1/4} N^{- 1/3}. \label{crg}
\end{equation}
Set $X := |\{x \in V(\Hat{G}^N) : |\mathcal{C} (x) | > \lambda^{- 1/12} N^{2/3}~\text{and}~\text{diam}_R (\mathcal{C} (x)) < \lambda^{-1} N^{1/3} \}|.$
Note that if $X > 0,$ then $X > \lambda^{- 1/12} N^{2/3}.$ So, by the Chebyshev inequality and (\ref{crg}), we have
\begin{align*}
&\textbf{P} (\exists x \in V(\Hat{G}^N), |\mathcal{C} (x) | > \lambda^{- 1/12} N^{2/3}~\text{and}~\text{diam}_R (\mathcal{C} (x)) < \lambda^{-1} N^{1/3}) \\
&\le \textbf{P} (X \ge \lambda^{- 1/12} N^{2/3} ) \\
&\le c_{11} \lambda^{-1/6}. ~\Box
\end{align*}

\textit{Proof of Proposition \ref{crprop}.} 
By Lemma \ref{crlem} (2) (4), (\ref{condition-of-type2-(1)-1}) and (\ref{condition-of-type2-(3)}) in Theorem \ref{type2} hold for $v(N) = N^{2/3}$ and $r(N) = N^{1/3}.$
So, we only need to check (\ref{condition-of-type2-(1)-2}) with $r(N) = N^{1/3}.$ 
The condition (\ref{condition-of-type2-(1)-2}) follows from Lemma \ref{crlem} and a minor modification of the proof of Theorem 3.1 of \cite{BDNP}.
To make the paper self-contained, we briefly recall the argument of \cite{BDNP}.
Fix $x \in V(\Hat{G}^N), 0 \le k \le k_0^N : = 2 \lfloor \log_2 \log N \rfloor,$ sufficiently large $\lambda > 0,$ and $N \ge \max \{ \lambda^4, a(\lambda) \}.$
By (\ref{ineqcr3}), we have a sequence $(r_{k, j}^N)_{j = 0}^{\lfloor 4 \lambda^2 2^k \rfloor}$ satisfying 
$r_{k, 0}^N = 0, \frac{(j - 1/2) N^{1/3}}{4 \lambda 2^k} \le r_{k, j}^N \le \frac{j N^{1/3}}{4 \lambda 2^k}$ and
$\textbf{E} |\partial B_p (x, r_{k, j}^N ; \Hat{G}^N ) | \le 16 \lambda^2 c_5 e^{A \lambda} 2^k$ for all $1 \le j \le \lfloor 4 \lambda^2 2^k \rfloor.$
We say that $y \in \partial B_p (x, r_{k, j}^N ; \Hat{G}^N)$ is $k$-good if $y$ and a vertex in $\partial B_p (x, r_{k, j+1}^N ; \Hat{G}^N)$ are connected by a path which does not
intersect $\partial B_p (x, r_{k, j}^N ; \Hat{G}^N)$ except $y.$
Set
$$A_{k}^N (x) := \begin{cases}
                           \displaystyle \bigcup_{j = 0}^{\lfloor 4 \lambda^2 2^k \rfloor} \{y \in \partial B_p (x, r_{k, j}^N ; \Hat{G}^N) : y~\text{is $k$-good} \}
 & \text{if} ~0 \le k \le k_0^N - 1, \\ 
                             \\
                            \mathcal{C} (x) & \text{if}~ k = k_0^N.
                            \end{cases}$$

We define $\ell_k^N(x) : = \frac{\text{diam}_R (\mathcal{C}(x))}{2^k}$ for $0 \le k \le k_0^N - 1$ and $\ell_{k_0^N}^N (x) := 0.$
Under the events that $\text{diam}(\mathcal{C}(x)) \le \lambda N^{1/3}$ and $\text{diam}_R (\mathcal{C} (x)) \ge \lambda^{-1} N^{1/3},$ 
the set of resistance balls $\{ B_{\text{eff}}^N (y, \ell_k^N (x)) : y \in A_k^N (x) \}$ is a covering for $\mathcal{C}(x)$
and $n_{\text{cov}} (\mathcal{C} (x), \ell_k^N (x)) \le |A_k^N (x)|$ for all $0 \le k \le k_0^N.$
By (\ref{ineqprobab}), we get for some $c_{12} > 0,$
\begin{align} 
&\textbf{P} \Big(\exists x \in V(\Hat{G}^N), 0\le \exists k \le k_0^N, |\mathcal{C} (x)| > \lambda^{-1/12} N^{2/3}, |A_k^N (x)| \ge e^{(|A| + 1) \lambda 2^{\frac{k}{2}}} \Big) \notag \\
&\le c_{12}^{-1} \lambda^{-c_{12}}. \label{ineq-covering}
\end{align}
Set $\ell_k^N := \frac{\text{diam}_R (\mathcal{C}^N)}{2^k}$ for $0 \le k \le k_0^N - 1$ and $\ell_{k_0^N}^N := 0.$
By (\ref{ineq-covering}) together with Lemma \ref{crlem} (2), (3) and (\ref{ineqcr2}), we have for some $c_{13}, c_{14} > 0,$
\begin{align*}
&\textbf{P} \Big(\sum_{k = 1}^{k_0^N} \sqrt{\ell_{k-1}^N \log \{ n_{\text{cov}} (\mathcal{C}^N, \ell_k^N ) \} } \ge c_{13} \lambda \sqrt{N^{1/3}}~ \Big) \\
&\le \textbf{P} \Big(|\mathcal{C}^N| > \lambda^{-1/12} N^{2/3}, \text{diam} (\mathcal{C}^N) < \lambda N^{1/3} \\
&~~~~~~~~~~~~~~~~~~~~~~~~~\text{and}~\sum_{k = 1}^{k_0^N} \sqrt{\ell_{k-1}^N \log \{ n_{\text{cov}} (\mathcal{C}^N, \ell_k^N ) \} } \ge c_{13}  \lambda \sqrt{N^{1/3}}~ \Big) \\
&+ \textbf{P} (|\mathcal{C}^N| \le \lambda^{-1/12} N^{2/3}) + \textbf{P} (\text{diam} (\mathcal{C}^N) \ge \lambda N^{1/3}) \\
&\le \textbf{P} (\exists x \in V(\Hat{G}^N), 
|\mathcal{C} (x) | > \lambda^{- 1/12} N^{2/3}~\text{and}~\text{diam}_R (\mathcal{C} (x)) < \lambda^{-1} N^{1/3} ) \\
&+ \textbf{P} \Big(\exists x \in V(\Hat{G}^N), 0 \le \exists k \le k_0^N, |\mathcal{C} (x)| > \lambda^{-1/12} N^{2/3}, |A_k^N (x)| \ge e^{(|A| + 1) \lambda 2^{\frac{k}{2}}} \Big) 
\\ 
&+ \textbf{P} (|\mathcal{C}^N| \le \lambda^{-1/12} N^{2/3}) + \textbf{P} (\text{diam} (\mathcal{C}^N) \ge \lambda N^{1/3}) \\
&\le c_{14}^{-1} \lambda^{-c_{14}}. \Box
\end{align*}

\subsection{The range of random walk in $\mathbb{Z}^{d}, d \ge 5$}
\label{sec:13}
~~~~Let $d \ge 5$. We write $(S_n)_{n \ge 0}$ to denote the simple random walk in $\mathbb{Z}^{d}$ started from 0 which is defined on a probability space with 
probability measure $\textbf{P}$.
Let $G^N$ be a graph with vertex set $V(G^N) := \{S_n : 0 \le n \le N \}$ and edge set $E(G^N) := \{ \{S_{n - 1}, S_n \} : 1 \le n \le N \}$. 
We prove the following proposition.
\begin{prop} \label{trprop}
There exist $c_1, c_2 > 0$ such that $\textbf{P}$-a.s., for sufficiently large $N \in \mathbb{N}$,
$$c_1 N^2 \le t_{\text{cov}}(G^N) \le c_2 N^2,$$
and $(t_{\text{cov}}(G^N))_{N \in \mathbb{N}}$ is Type 2.
\end{prop}
Let $(S_{- n} )_{n \ge 0}$ be an independent copy of $(S_n)_{n \ge 0}$ and set $S = (S_n)_{n \in \mathbb{Z} }$.
Let $\mathcal{T}$ be the set of cut-times, that is, $\mathcal{T} := \{n : S_{(- \infty, n] } \cap S_{[n + 1, \infty) } = \emptyset \}.$ 
We can write $\mathcal{T} \cap (0, \infty) = \{T_n : n \in \mathbb{N} \}$, where $0 < T_1 < T_2 < \dotsc $. Set cut-points $C_n := S_{T_n}$. 
We use the following fact.
\begin{lem} \label{RRW} (\cite{Cr}, Lemma 2.2 (see also \cite{CHK}, (5.6)) )  
$$\lim_{n \to \infty} \frac{T_n}{n} = \tau(d) := \mathbf{E}(T_1 | 0 \in \mathcal{T} ) \in [1, \infty), ~~\textbf{P}-a.s.$$ 
\end{lem}
\textit{Proof of Proposition \ref{trprop}.} 
We check almost-sure versions of (\ref{condition-of-type2-(1)-1}), (\ref{condition-of-type2-(1)-2}) and
(\ref{condition-of-type2-(3)}) in Theorem \ref{type2} with $v(N) = r(N) = N.$
For $N \in \mathbb{N}$, there exists $M = M(N) \in \mathbb{N}$ such that $T_M \le N < T_{M + 1}$. 
Because $d_{G^N}(0, C_M) \ge M,$ we have that $|E(G^N)| \ge M, ~~\textbf{P}$-a.s.
By Lemma \ref{RRW}, there exist $c_3, c_4 > 0$ such that 
$c_3 N \le M \le c_4 N,$ for sufficiently large $N \in \mathbb{N}$, $\textbf{P}$-a.s. 
So, $\textbf{P}$-a.s., for sufficiently large $N \in \mathbb{N}$, 
$$|E(G^N)| \ge c_3 N.$$
Every path from $0$ to $C_M$ must pass edges $\{S_{T_n}, S_{T_n + 1} \}_{1 \le n \le M - 1}$.
So, by Lemma \ref{Nash-Williams}, there exists $c_5 > 0$ such that $\textbf{P}$-a.s., for sufficiently large $N \in \mathbb{N}$, 
\begin{equation}
\text{diam}_R(G^N) \ge R_{\text{eff}}^N(0, C_M) \ge M - 1 \ge c_5 N. \label{ineqrr}
\end{equation}
By definition, 
$$|E(G^N)| \le N, ~\text{and}~ \text{diam}_R(G^N) \le \text{diam}(G^N) \le N, ~ \textbf{P}-a.s.$$   
Fix $1 \le k \le \lfloor \log_2 \log (c_5 N) \rfloor.$ We define $A_k^N$ as follows:
$$A_k^N := \begin{cases} \{S_{\lfloor j \frac{c_5 N}{2^{k + 1}} \rfloor} : 0 \le j \le \lfloor \frac{2^{k + 1}}{c_5} \rfloor \}, & \text{if}~1 \le k \le \lfloor \log_2 \log (c_5 N) \rfloor - 1, \\
                                  \{ S_j  : 0 \le j \le N \} & \text{otherwise}.
               \end{cases}$$
It is not hard to check that $V(G^N) \subset \bigcup_{u \in A_k^N } B^N (u, \frac{c_5 N}{2^k} ),$ 
where $B^N (u, r ) = \{v \in V(G^N) : d_{G^N} (u, v) \le r \}.$
Set $k_0^N = \lfloor \log_2 \log (c_5 N) \rfloor.$
By (\ref{ineqrr}),
we have that $\textbf{P}$-a.s., for sufficiently large $N \in \mathbb{N}$,
$$V(G^N) \subset \bigcup_{u \in A_k^N } B_{\text{eff} }^N (u, \ell_k^N ),$$
where $\ell_k^N = \frac{\text{diam}_R (G^N)}{2^k}$ for $1 \le k \le k_0^N - 1$ and $\ell_k^N = 0$ otherwise.
Because $n_{\text{cov}} (G^N, \ell_k^N) \le |A_k^N| \le \lfloor \frac{2^{k + 1}}{c_5} \rfloor + 1 \le c_6 2^k$ for some $c_6 > 0$ and all $k < k_0^N,$
we have $\textbf{P}$-a.s., for sufficiently large $N \in \mathbb{N},$
$$\sum_{k=1}^{k_0^N} \sqrt{\ell_{k-1}^N \log \{n_{\text{cov}} (G^N, \ell_k^N) \}} \le c_7 \sqrt{N} ~\text{for some}~c_7 > 0.$$
By Remark \ref{type2rems} (3), we complete the proof. $\Box$

\subsection{Sierpinski gasket graphs}
\label{sec:14}
~~~~Let $p_1, p_2, p_3$ be vertices of an equilateral triangle in $\mathbb{R}^2$. We define three contraction maps $\psi_i : {\mathbb{R} }^2 \to {\mathbb{R} }^2, i = 1, 2, 3$ as follows:

$$\psi_i (x) = p_i + \frac{x - p_i }{2}, ~~~i = 1, 2, 3, x \in {\mathbb{R} }^2.$$

$G^N$ is a graph with the following vertex and edge sets:

$V(G^N) := \displaystyle \bigcup_{i_1 \dotsc i_N =1 }^{3} \psi_{i_1 \dotsc i_N } (V_0), $ \\

$E(G^N) := \{\{\psi_{i_1 \dotsc i_N } (x), \psi_{i_1 \dotsc i_N } (y)  \} : x, y \in V_0, x \neq y, i_1, \dotsc, i_N \in \{1, 2, 3 \} \}$, \\

where $V_0 := \{p_1, p_2, p_3 \}$ and $\psi_{i_1 \dotsc i_N } := \psi_{i_1} \circ \dotsc \circ \psi_{i_N}.$ \\
Random weights $(\mu_{xy}^N )_{\{x, y \} \in E(G^N) }$ are i.i.d. random variables with a common distribution which is supported on $[c_1, c_2]$, where $0< c_1 \le c_2 < \infty.$
We will establish the following estimate of the cover time for $G^N$:
\begin{prop} \label{sgprop} There exist $c_3, c_4 > 0$ such that for all $N \in \mathbb{N}$, $\textbf{P}$-a.s., 
$$c_3 5^N \le t_{\text{cov} } (G^N) \le c_4 5^N,$$
and $(t_{\text{cov}} (G^N))_{N \in \mathbb{N}}$ is Type 2.
\end{prop}

To prove this proposition, we prepare some notations.
For $i_1, \dotsc, i_n \in \{1, 2, 3 \}$ and $n \le N,$ let $G_{i_1 \dotsc i_n}^N$ be the induced graphs with vertex set
$V(G_{i_1 \dotsc i_n}^N)$ which is the intersection of $V(G^N)$ and an equilateral triangle with vertices
$\psi_{i_1 \dotsc i_n}(p_i), i = 1, 2, 3.$

We use the following lemma. The resistance estimate is obtained, for example, from arguments in section 7 of \cite{Ba} or section 1.3 of \cite{St}.
\begin{lem} \label{sglem} There exist $c_5, c_6 > 0$ such that for all $N \in \mathbb{N},$
$$c_5 3^N \le |\mu(G^N)| \le c_6 3^N,~~~c_5 \left(\frac{5}{3} \right)^N  \le \text{diam}_R(G^N) \le c_6 \left(\frac{5}{3} \right)^N ~\textbf{P}\text{-a.s.}$$
\end{lem}              
                           
\textit{Proof of Proposition \ref{sgprop}.}
By Lemma \ref{sglem}, almost-sure versions of (\ref{condition-of-type2-(1)-1}) and (\ref{condition-of-type2-(3)}) hold for $v(N) = 3^N$ and $r(N) = (\frac{5}{3})^N.$
We only need to check an almost-sure version of (\ref{condition-of-type2-(1)-2}) with $r(N) = (\frac{5}{3})^N.$ \\
Set $\ell_k^N = c_6 (\frac{5}{3})^{N - k}$ for $0 \le k < N$ and $\ell_k^N = 0$ otherwise. 
Let $x_{i_1, \cdots, i_k}^N$ be a fixed vertex in $V(G_{i_1 \cdots i_k}^N).$
By Lemma \ref{sglem}, $\{B_{\text{eff}}^N (x_{i_1 \cdots i_k}^N, \ell_k^N ) : i_1, \cdots, i_k \in \{1, 2, 3 \} \}$ is a covering for $G^N$
$\textbf{P}\text{-a.s.}$ In particular, we get
$$n_{\text{cov}} (G^N, \ell_k^N ) \le 3^k ~\textbf{P}\text{-a.s.}$$
Therefore, we have for some $c_7 > 0$ and all $N \in \mathbb{N},$
$$\sum_{k = 1}^{N} \sqrt{\ell_{k-1}^N \log \{ n_{\text{cov}} (G^N, \ell_k^N) \}} \le c_7 \sqrt{\Big (\frac{5}{3} \Big )^N} ~\textbf{P}\text{-a.s.}$$
By Remark \ref{type2rems} (3), we complete the proof. $\Box$ 
\begin{rem}
It will be possible to estimate cover times for Sierpinski gasket graphs in higher dimensions and nested fractals by applying arguments similar to the above proof. 
\end{rem}
{\bf Acknowledgements.} \\
We would like to thank Professor Takashi Kumagai for fruitful discussions and careful reading of an early version of this paper. 
We would like to thank Dr. Ryoki Fukushima for helpful comments and suggestion that the Seneta-Heyde theorem (Lemma \ref{scgwlem} (2)) is useful 
for the proof of Proposition \ref{scgwprop}.

\end{document}